\documentclass[12pt]{amsart}
\date{6 January 2019}
\usepackage[T1]{fontenc}
\usepackage{latexsym,amsmath,amscd,amssymb,lmodern,amsfonts}
\usepackage{lscape}
\usepackage{stmaryrd}
\usepackage{array}
\usepackage{pdfsync}
\usepackage[all]{xy}
\usepackage[english]{babel}

\usepackage[mathscr]{eucal} 
\usepackage{mathrsfs}
\usepackage{graphicx}
\usepackage{calligra}

\usepackage{amsthm} 
\usepackage{xfrac}  
\usepackage{multirow} 
\usepackage{setspace}

\usepackage{hyperref}
\hypersetup{colorlinks=true,
    linkcolor=black,          
    citecolor=black,        
    filecolor=black,      
    urlcolor=cyan           
}

\def\addsymbol #1: #2#3{$#1$ \> \parbox{5in}{#2 \dotfill \pageref{#3}}\\} 
\def\newnot#1{\label{#1}}

\newcommand{\st}{\;|\;}
\newcommand{\xra}{\xrightarrow}

\newcommand{\ts}{\times}

\newcommand{\ot}{\otimes}

\renewcommand{\H}{\textrm{H}}

\newcommand{\Z}{{\mathbb{Z}}}
\newcommand{\Q}{{\mathbb{Q}}}
\newcommand{\R}{{\mathbb{R}}}
\newcommand{\C}{{\mathbb{C}}}
\newcommand{\Oc}{{\mathbb{O}}}

\newcommand{\cD}{\mathcal{D}}

\newcommand{\cM}{\mathcal{M}}
\newcommand{\cO}{\mathcal{O}}
\newcommand{\cR}{\mathcal{R}}

\newcommand{\chS}{\check{S}}

\newcommand{\fe}{\mathfrak{e}}

\newcommand{\fl}{\mathfrak{l}}
\newcommand{\fm}{\mathfrak{m}}

\newcommand{\fp}{\mathfrak{p}}

\newcommand{\fu}{\mathfrak{u}}

\newcommand{\fso}{\mathfrak{so}}

\newcommand{\liea}{\mathfrak{a}}

\newcommand{\liec}{\mathfrak{c}}
\newcommand{\lied}{\mathfrak{d}}
\newcommand{\liem}{\mathfrak{m}}

\newcommand{\liep}{\mathfrak{p}}

\newcommand{\lies}{\mathfrak{s}}
\newcommand{\liet}{\mathfrak{t}}
\newcommand{\liemp}{\mathfrak{m}_+}

\newcommand{\liemc}{\mathfrak{m}^{\mathbb{C}}}
\newcommand{\lieh}{\mathfrak{h}}
\newcommand{\liehc}{\mathfrak{h}^{\mathbb{C}}}
\newcommand{\lieg}{\mathfrak{g}}
\newcommand{\liel}{\mathfrak{l}}
\newcommand{\liez}{\mathfrak{z}}
\newcommand{\liesu}{\mathfrak{su}}
\newcommand{\lieso}{\mathfrak{so}}
\newcommand{\liesl}{\mathfrak{sl}}
\newcommand{\lieu}{\mathfrak{u}}

\newcommand{\tlieg}{\widehat{\lieg}}

\newcommand{\al}{\alpha}
\newcommand{\be}{\beta}

\newcommand{\ga}{\gamma}

\newcommand{\De}{\Delta}
\newcommand{\Ga}{\Gamma}

\newcommand{\Sh}{\check{S}}

\newcommand{\PSL}{\mathrm{PSL}}

\newcommand{\PSU}{\mathrm{PSU}}
\newcommand{\SU}{\mathrm{SU}}
\newcommand{\U}{\mathrm{U}}
\newcommand{\GL}{\mathrm{GL}}
\newcommand{\SL}{\mathrm{SL}}
\newcommand{\SO}{\mathrm{SO}}
\newcommand{\Sp}{\mathrm{Sp}}
\newcommand{\SSS}{\mathrm{S}}
\newcommand{\Spin}{\mathrm{Spin}}
\newcommand{\E}{\mathrm{E}}
\newcommand{\F}{\mathrm{F}}
\newcommand{\G}{\mathrm{G}}

\newcommand{\OO}{\mathrm{O}}
\newcommand{\Or}{\mathrm{O}}

\DeclareMathOperator{\ad}{ad}
\DeclareMathOperator{\Ad}{Ad}

\DeclareMathOperator{\rk}{rk}
\DeclareMathOperator{\rank}{rank}
\DeclareMathOperator{\reg}{reg}

\DeclareMathOperator{\Hom}{Hom}

\DeclareMathOperator{\End}{End}

\DeclareMathOperator{\Id}{Id}

\DeclareMathOperator{\Herm}{Herm}

\DeclareMathOperator{\Aut}{Aut}

\DeclareMathOperator{\Mat}{Mat}

\DeclareMathOperator{\Sym}{Sym}
\DeclareMathOperator{\Skew}{Skew}

\DeclareMathOperator{\Tr}{Tr}

\DeclareMathOperator{\Lie}{Lie}
\DeclareMathOperator{\Split}{split}
\DeclareMathOperator{\opp}{opp}

\newcommand{\lra}{\longrightarrow}

\newcommand{\rmin}{\mathrm{top}}

\newcommand{\noi}{\noindent}

\newcommand{\abs}[1]{\lvert#1\rvert}

\newcommand{\lie}[1]{\mathfrak{#1}}

\newcommand{\sll}{\mathfrak{sl}}
\newcommand{\projects}{\twoheadrightarrow}

\newcommand{\mtrx}[1]{\left (\begin{matrix}#1\end{matrix}\right)}
\newcommand{\smtrx}[1]{\left (\begin{smallmatrix}#1\end{smallmatrix}\right)}

\newtheorem{theorem}{Theorem}[section]
\newtheorem{corollary}[theorem]{Corollary}
\newtheorem{lemma}[theorem]{Lemma}
\newtheorem{proposition}[theorem]{Proposition}

\numberwithin{equation}{section}

\newtheorem*{theorem*}{Theorem}

\newtheorem{definition}[theorem]{Definition}

\theoremstyle{remark}
\newtheorem{remark}[theorem]{Remark}

\begin{document}

\title[Higgs bundles and higher Teichm\"uller spaces]
{Higgs bundles and higher Teichm\"uller spaces}

\author[Oscar Garc{\'\i}a-Prada]{Oscar Garc{\'\i}a-Prada}
\address{Instituto de Ciencias Matem\'aticas \\
  CSIC-UAM-UC3M-UCM \\ Nicol\'as Cabrera, 13--15 \\ 28049 Madrid \\ Spain}
\email{oscar.garcia-prada@icmat.es}

\thanks{
  The  author is partially supported by the Spanish MINECO under 
the ICMAT Severo Ochoa grant No. SEV-2011-0087, and under grant 
No. MTM2013-43963-P.  
}


\begin{abstract}
This paper is a survey on  the role of Higgs bundle theory in the study of higher Teichm\"uller spaces. Recall that the Teichm\"uller space 
of a compact surface can be identified with a certain connected component of the moduli space of representations of the fundamental group 
of the surface into $\PSL(2,\R)$. Higher Teichm\"uller spaces correspond to special components  of the moduli space of representations when 
one  replaces $\PSL(2,\R)$ by a real non-compact semisimple Lie group of higher rank. Examples of these spaces are provided by the Hitchin 
components for  split real groups, and maximal Toledo invariant components for groups of Hermitian type.  More recently, the existence of 
such components has been proved for $\SO(p,q)$, in agreement with  the conjecture of  Guichard and Wienhard relating the  existence of 
higher Teichm\"uller spaces to a certain  notion of positivity on a Lie group that they have  introduced. We review these three different 
situations, and end up  explaining briefly the conjectural general picture from the point of view of Higgs bundle theory. 
\end{abstract}
\maketitle

\section{Introduction}

Given a compact surface $S$ of genus $g\geq 2$,  consider the moduli space of 
representations of the fundamental group of $S$ in $\SL(2,\R)$.  Attached to such 
a representation there is an integer invariant $d$, which by Milnor \cite{milnor} 
satisfies the  bound $|d|\leq g-1$. In the maximal case $d=g-1$ (or $d=1-g$) the moduli space
has $2^{2g}$ connected components, all of which are homeomorphic to $\R^{6g-6}$ and can
be identified to the Teichm\"uller space of $S$. All the representations in these components are
Fuchsian, that is, discrete and faithful. When passing to the adjoint group 
$\PSL(2,\R)$ all these components get identified.

Higher Teichm\"uller spaces appear when one replaces $\SL(2,\R)$ by certain real Lie groups of
higher rank. These spaces have been studied from several points of view. 
In this paper we will focus on the role played by the theory of Higgs bundles.  
This powerful theory, introduced by Hitchin \cite{hitchin1987}, is very useful in 
identifying these spaces and studying their topology, in particular counting the number of connected 
components. To do this, one chooses a complex 
structure on $S$, making it into a compact Riemann surface $X$. Given a semisimple real Lie group $G$, 
with a choice of a maximal compact subgroup $H\subset G$, one has a Cartan decomposition of the Lie algebra
of $G$, given by $\lieg=\lieh \oplus \liem$, where $\lieh$ is the Lie algebra of $H$ and $\liem$ is the 
orthogonal complement of $\lieh$ with respect to the Killing form of $\lieg$. A $G$-Higgs bundle over 
$X$ is a pair $(E,\varphi)$ consisting of a holomorphic  $H^\C$-bundle $E$ over $X$ and  a holomorphic 
section $\varphi$ of the bundle $E(\liem^\C)\otimes K$, where $E(\liem^\C)$ is the vector bundle associated to $E$ via the
isotropy representation $H^\C\to \GL(\liem^\C)$, and $K$ is the canonical line bundle of $X$ --- the holomorphic 
cotangent bundle of $X$. The non-abelian Hodge correspondence proved  by Hitchin \cite{hitchin1987}, Donaldson \cite{donaldson},
Simpson \cite{Sim88,Sim92}, Corlette \cite{corlette} and others (see \cite{GGM09}) establishes a homeomorphism between the moduli 
space  $\cM(G)$ of polystable $G$-Higgs bundles over $X$ and the moduli space  $\cR(G)$
of reductive representations of the fundamental group of $S$ in $G$.

In  \cite{hitchin92} Hitchin constructed components of the moduli space of representations $\cR(G)$ when $G$ is the split 
real form of a complex simple Lie group $G^\C$ which, like Teichm\"uller space, are homeomorphic to a real vector space, in
this case of dimension $\dim G(2g-2) $. Originally referred by him as generalized Teichm\"uller components, they are now  
called Hitchin components. Such  components appear as  sections 
of the Hitchin fibration for the moduli space 
of $G^\C$-Higgs bundles. This fibration is defined by the  map $h:\cM(G^\C)\to B(G^\C)$ sending a $G^\C$-Higgs bundle 
$(E,\varphi)$ to $(p_1(\varphi),\ldots,p_r(\varphi))$, where 
$p_1,\ldots,p_r$ is   a basis of the ring  of
$G^\C$-invariant homogeneous polynomials of $\lieg^\C$. Hence  $B(G^\C)= \oplus_{i=1}^rH^0(X,K^{d_i})$, where $d_i$ is the degree of $p_i$. 
These degrees are indepedent of the basis and depend only on the group. The Hitchin components can then be identified with the vector space
 $B(G^\C)$ which has half the dimension of $\cM(G^\C)$. An alternative (equivalent)  way of definining the Hitchin components is by directly considering the Hitchin map for the Lie group $G$. In fact this exists for any semisimple Lie group $G$. The Hitchin map in this case 
$h:\cM(G)\to B(G)$ is obtained by applying a basis of the ring of $H^\C$-invariant polynomials of $\liem^\C$. One can construct in this general situation a section of $h$, but only when $G$ is split does this define a component of $\cM(G)$
--- when $G$ is split $B(G)=B(G^\C)$ (see \cite{GPR18}). 

The Hitchin components  consist entirely of representations of $\pi_1(S)$ in $G$ that can be deformed to a 
representation that factors through a Fuchsian representation  $\pi_1(S)\to \SL(2,\R)$, where $\SL(2,\R)\to G$ is the 
irreducible representation of  $\SL(2,\R)$ in $G$, which  always exists when $G$ is a split real form.  One may need to consider 
$\PSL(2,\R)$ if the group $G$  is of adjoint type. In this case, like for $\PSL(2,\R)$,  the Hitchin component is unique.
It follows from the work of Labourie  \cite{labourie}  and Fock and Goncharov \cite{fock-goncharov} that, similarly to the 
Teichm\"uller component for $\PSL(2,\R)$, all the representations in a Hitchin component are discrete and faithful.
In his paper  \cite{hitchin92} Hitchin already posed the question on the relation of his generalized Teichm\"uller components 
with geometric structures. A  geometric interpretation of these components for $G=\SL(3,\R)$ has been given by Choi--Goldman 
\cite{ChG93} and more recently for more general split groups by Labourie \cite{labourie} Fock--Goncharov 
\cite{fock-goncharov}, and Guichard--Wienhard \cite{guichard-wienhard:2008}.

The group $\SL(2,\R)$ is of course a split real form of $\SL(2,\C)$ but it is also a group of Hermitian type. For such a group $G$, 
the symmetric space $G/H$ is K\"ahler.  The symmetric space $\SL(2,\R)/\SO(2)$  can be identified with the Poincar\'e upper-half plane. Another direction in which Teichm\"uller space can be generalized is 
by considering a higher rank non-compact simple Lie group $G$ of Hermitian type 
with finite centre. The Lie algebra of a group of this type is
one of the following: $\mathfrak{su}(p,q)$, $\fso^*(2n)$, 
$\mathfrak{sp}(2n,\R)$, $\fso(2,n)$, $\fe_6^{-14}$ and $\fe_7^{-25}$
(we are using Helgason's notation \cite{helgason}).
In this situation 
the centre $\liez$ of $\lieh$ is isomorphic to $\R$, and  
the adjoint action of a special element $J\in \liez$ defines an almost complex 
structure on $\liem=T_o(G/H)$, where $o\in G/H$ corresponds to the coset $H$, 
making the  symmetric space $G/H$ into a K\"ahler manifold. 
The almost complex structure $\ad(J)$ gives a decomposition 
$\liem^\C=\liem^+ + \liem^-$ in $\pm i$-eigenspaces, which is $H^\C$-invariant.
An immediate consequence of  this decomposition  for a $G$-Higgs bundle 
$(E,\varphi)$ is that it  gives a 
bundle decomposition $E(\liem^\C)=E(\liem^+) \oplus E(\liem^-)$ 
and hence the Higgs field decomposes as $\varphi=(\varphi^+,\varphi^-)$,
where $\varphi^+\in H^0(X,E(\liem^+)\otimes K)$ and 
$\varphi^-\in H^0(X,E(\liem^-)\otimes K)$.

Another  important ingredient of $G$-Higgs bundles for this class
of groups is the Toledo invariant. This is a topological 
invariant attached to a $G$-Higgs bundle $(E,\varphi)$, in fact attached to 
$E$. It is defined   by considering a special character $\chi_T$  of 
$\lieh^\C$ called the Toledo character. 
If this lifts to a character $\tilde{\chi}_T$
of $H^\C$, we consider the associated line bundle $E(\tilde{\chi}_T)$, and 
define the Toledo invariant $\tau$ of $(E,\varphi)$ as
$$
\tau=\tau(E):=\deg(E(\tilde{\chi}_T)).
$$  
Otherwise one can show that there is a 
rational number $q_T$ such $q_T\chi_T$ lifts to a character
$\tilde{\chi}_T$, and one can  define 
$\tau=\tau(E):=\frac{1}{q_T}\deg(E(\tilde{\chi}_T))$.  A crucial fact is given
by the Milnor--Wood inequality: The moduli space $\cM(G)$ is empty unless 
$$
|\tau|\leq \rk(G/H)(2g-2).
$$
This was proved on a case by case  basis for the classical groups 
\cite{hitchin1987,gothen,BGG03,BGG06,BGG15,GGM13} and in general in \cite{BGR17}. 
The Toledo invariant of $(E,\varphi)$ coincides with the classical Toledo invariant 
of a representation of the fundamental group in $G$, which is defined 
by integrating over $X$ the pull-back of the K\"ahler form of the Bergman metric
of $G/H$, via a $\rho$-equivariant map $f:\widetilde{X}\to G/H$ determined by
$\rho$. The map $f$ can be taken to be harmonic if the representation is
reductive, and hence corresponding to a polystable $G$-Higgs bundle.
In the context of representations the inequality $|\tau|\leq \rk(G/H)(2g-2)$,
goes back to Milnor \cite{milnor}, who studies the case $G=\SL(2,\R)$, as mentioned above,  
and was proved  in various cases in  \cite{wood,dupont,DT87,CO03}, and in general 
by Burger--Iozzi--Wienhard \cite{BIW10}. 

The higher Teichm\"uller spaces appear when the Toledo invariant is maximal, i.e., when 
$|\tau|= \rk(G/H)(2g-2)$, and the symmetric space $G/H$ is of tube type.  A geometric characterization
of the tube type condition  is
given in terms of the Shilov boundary of the Harish-Chandra bounded symmetric domain
realization of $G/H$. We say that $G/H$ is of tube type if the Shilov boundary
of the corresponding bounded symmetric domain is a symmetric space $H/H'$ (of
compact type). 
If this is the case, this domain is biholomorphic to a `tube' over the
symmetric  cone $\Omega:=H^*/H'$ --- the 
non-compact symmetric space dual of the Shilov boundary. The simplest example
of this is given by $G=\SL(2,\R)\cong \SU(1,1)$, where the bounded symmetric domain 
is the Poincar\'e disc and  the tube realization is  
the Poincar\'e upper half plane (these are of course related by the classical Cayley transform).
 Sometimes we will say that the group $G$ is of tube type
if the symmetric space $G/H$ is of tube type. The groups of tube type have 
Lie algebra given by $\mathfrak{su}(p,q)$ with $p=q$, $\fso^*(2n)$ with $n$
even, $\mathfrak{sp}(2n,\R)$, $\fso(2,n)$  
and $\fe_7^{-25}$. While those of non-tube type are $\mathfrak{su}(p,q)$ with  
$p\neq q$, $\fso^*(2n)$ with $n$ odd, and $\fe_6^{-14}$.

The subvariety $\cR_{\max}(G)\subset \cR(G)$ of maximal representations  
(and corresponding  subvariety  $\cM_{\max}(G)\subset \cM(G)$) in the tube case
has special significance since, as proved in \cite{BIW10},  it consists entirely of
discrete and faithful representations, like the Hitchin components in the split case.  
As in this case, all the representations in $\cR_{\max}(G)$ are  Anosov in the sense 
of Labourie \cite{labourie,burger-iozzi-labourie-wienhard:2005}, and are related to 
geometric structures of various kinds \cite{GW12}. 

From the Higgs bundle point of view the  tube-type condition plays a fundamental role
in constructing a bijective correspondence --- called Cayley correspondence --- between 
the moduli space of maximal Higgs bundles 
$\cM_{\max}(G)$ and the moduli space of $K^2$-twisted $H^*$-Higgs bundles over $X$, where 
$H^*\subset H^\C$ is the non-compact dual of $H$, described above. These are defined like 
Higgs bundles except that the twisting is by $K^2$ instead of $K$. An immediate consequence 
of this correspondence is the  existence of new topological invariants for maximal Higgs 
bundles   which are hidden a priori, and which play a  crucial role in determining the number 
of connected components of $\cM_{\max}(G)$.

Up to a finite cover, there is only one group which is both split and of 
Hermitian type, and in fact of tube type. This is the symplectic group 
$\Sp(2n,\R)$. The tube realization is the Siegel upper half space. The Hitchin 
representations are actually maximal but, except for the case $n=1$, there are
maximal components that are not Hitchin components (see \cite{gothen,GGM13}).

Maximal Higgs bundles, and hence maximal representations, in the non-tube case 
present also very interesting phenomena.  It turns out that in this case the dimension 
of the moduli space of maximal $G$-Higgs bundles is smaller than expected, since up to 
a compact factor every maximal representation reduces to a  representation in a maximal 
subgroup of tube type. This rigidity phenomenon is  very rare in the context of  surface 
groups (see  \cite{toledo,hernandez,BGG03,BGG06, BIW10, BGG15,BGR17}), 
and it is conjectured to happen only for Hermitian groups of non-tube type \cite{KP14}.

For a while, split real groups and Hermitian groups were the only cases providing examples of 
higher  Teichm\"uller spaces.  However, gradually there was evidence that there may be other 
groups. Groups under scrutiny  were the special orthogonal groups of signature $(p,q)$. 
Of course, for certain values of $p$ and $q$  the group  is split or Hermitian, but for 
most values is neither.  One such evidence came from the  work of Guichard and Wienhard 
\cite{GW16}. They defined a notion of positivity for certain 
Lie groups and for surface group representations into those groups. The notion of positivity 
for a representation, which refines  Labourie's Anosov condition, is open and they conjecture 
that is also closed, hence detecting connected components of the moduli
space of representations. They showed that apart from the split real forms and the real forms of 
Hermitian type, the only other non-exceptional groups which allow positive representations are the 
groups locally isomorphic  to $\SO(p,q)$ for $1<p<q$. Indeed,  Collier \cite{Col17} has proved 
existence of special components containing positive representations for $\SO(n,n+1)$. Although this
group is split, the  Collier components are different to the  Hitchin ones. More recently the existence 
of special components  for $\SO(p,q)$ with  $1<p<q$ containing positive representations has been shown in \cite{A-et-al1,A-et-al2}.

After discribing the basic ingredients in the theory of Higgs bundles and the three different situations
where higher Teichm\"uller spaces emerge (split real forms, real forms of Hermitian type and $\SO(p,q)$), 
we finish the paper with some comments on the conjectural general picture. This seems to relate
positive representations, and hence higher Teichm\"uller spaces,  to components of the moduli space for which
there exists a Cayley correspondence  that generalizes the one for Hermitian groups of tube type and the one 
for the special components in the moduli space of  $\SO(p,q)$-Higgs bundles.
Indeed there is also a Cayley correspondence for the  Hitchin components.  This seems to be the way in which 
positivity is manifesting from the point of view of Higgs bundle theory. A full study of this is under investigation
in \cite{BCGGO}, including the exceptional groups for which there is a notion of positivity. As shown in
\cite{GW16} these  are real forms of $F_4$, $E_6$, $E_7$ and $E_8$ whose restricted root system is of type $F_4$.

As mentioned above, we have focused here on the point of view provided by the 
theory of Higgs bundles, and offering just a very partial picture. For other 
approaches,  and additional references, 
one may  see the very nice survey papers \cite{BIW14,wienhard}. We have not touched at all, the case of non-compact 
surfaces for which there is also a very rich  higher Teichm\"uller theory. This situation is studied for example in  \cite{fock-goncharov,BIW10}.
The Higgs bundle approach for the case of surfaces with punctures is very little developped. This involves the theory
of parabolic Higgs bundles. In \cite{BGM15} a non-abelian Hodge correspondence is established between parabolic $G$-Higgs bundles  
and representations of the fundamental group of the punctured surface with fixed conjugacy classes around the puntures, where 
$G$ is any real form of a complex reductive Lie group. This generalizes the work of Simpson \cite{Sim90} when $G=\GL(n,\C)$.
Although the split  and the Hermitian cases are briefly described  in  
\cite{BGM15}, much work has to be done along the lines 
of the compact case to study the higher Teichm\"uller spaces.

{\bf Acknowledgements}. The author would like to thank Athanase Papadopoulos for his kind invitation to contribute to this volume 
with this paper, and his patience in the delivery. He also wishes to thank his collaborators in the works on which this
paper is based, and Nigel Hitchin for inspiration, suggestions and many discussions over 
the years.

\section{Higgs bundles}
\label{chap:higgs-bundles}
\label{sec:g-higgs-defs}
\label{sec:moduli-spaces}
\label{sec:G-Higgs-herm}

\subsection{Basic definitions}

Following Knapp \cite[VII.2]{knapp}, we define a  \textbf{real reductive 
group} is a $4$-tuple $(G,H,\theta,B)$ where
\begin{enumerate}
\item $G$ is a real Lie group with reductive Lie algebra $\lieg$.
\item $H<G$ is a maximal compact subgroup.
\item $\theta$ is a Lie algebra involution of $\lieg$ inducing an eigenspace decomposition
$$\lieg=\lieh\oplus\liem$$
where $\lieh=\mathrm{Lie}(H)$ is the $(+1)$-eigenspace for the action of $\theta$, and $\liem$ is the 
$(-1)$-eigenspace.
\item\label{bilinear form} $B$ is a $\theta$- and $\Ad(G)$-invariant non-degenerate bilinear form,
with
respect to which $\lieh\perp_{B}\liem$
and $B$ is negative definite  on $\lieh$ and positive definite on $\liem$.
\item\label{diffeo} The multiplication map $H\times \exp(\liem)\to G$ is a 
diffeomorphism.

\item  $G$ acts by  inner automorphisms on the complexification $\lieg^\C$ of its Lie algebra via the adjoint representation. 
\end{enumerate}

In this section we fix a real reductive Lie group $(G,H,\theta,B)$.
We have $[\liem,\liem]\subset \lieh$, $[\liem,\lieh]\subset\lieh$.
Complexifying the isotropy representation $H\to \Aut(\liem)$, we obtain the
representation $\Ad: H^\C \to \Aut(\liem^\C)$.
When $G$ is semisimple we take $B$ to be the Killing form.
In this case $B$ and a choice of a maximal compact subgroup $H$ 
determine a Cartan decomposition.

Let $X$ be  a compact  Riemann surface of genus $g$. 

A {\bf $G$-Higgs bundle} on  $X$ consists of a holomorphic principal
$H^\C$-bundle  $E$ together with a holomorphic section $\varphi\in
H^0(X,E(\liem^\C)\ot K)$,
where  $E(\liem^\C)$ is the associated vector bundle with fibre 
$\liem^\C$ via the complexified isotropy representation, and $K$ is 
the canonical line bundle of $X$.

If $G$ is compact,  $H=G$ and $\liem=0$. A $G$-Higgs bundle is hence 
simply  a holomorphic principal  $G^\C$-bundle.
If  $G=H^\C$, where  $H$ is 
a maximal compact subgroup of $G$, $\liem=i\lieh$, and a 
$G$-Higgs  bundle is a principal $G$-bundle together with a section 
$\varphi\in H^0(X, E(\lieg) \ot K)$, where 
$E(\lieg)$ is the adjoint bundle. This is the original definition for 
complex Lie groups given by Hitchin in \cite{hitchin:duke}.

Let $G'$ be a reductive subgroup of $G$. A maximal compact subgroup of $G'$ 
is given by $H'=H\cap G'$ and we can take a compatible Cartan decomposition, 
in the sense that $\lieh'\subset \lieh$ and $\liem'\subset \liem$. Moreover, 
the isotropy representation of $H'$ is the restriction of the isotropy
representation $\Ad$ of $H$. We say that the structure group of a $G$-Higgs 
bundle 
$(E,\varphi)$ reduces to $G'$ when there is a reduction of the structure group 
of  the underlying $H^\C$-bundle to $H'^\C$, given by a subbundle $E_\sigma$, 
and  the Higgs field $\varphi\in H^0(X,E(\liem^\C)\ot K)$ belongs to 
$H^0(X,E_\sigma(\liem'^\C)\ot K)$.

\subsection{Stability of $G$-Higgs bundles}
\label{subsec:stability-G-Higgs-bundles}

There is a  notion of stability for $G$-Higgs
bundles (for details see \cite{GGM09}). To explain 
this  we  consider the   parabolic subgroups of $H^\C$ 
defined for $s\in i\lieh$ as
$$
P_s =\{g\in H^\C\ :\ e^{ts}ge^{-ts}\text{ is bounded as }t\to\infty\}.
$$
When $H^\C$ is connected every parabolic subgroup is conjugate to one of the 
form $P_s$ for some $s\in i\lieh$, but this is not the case necessarily when 
$H^\C$ is non-connected.
A Levi subgroup of $P_s$ is given by $L_s =\{g\in H^\C\ :\ \Ad(g)(s)=s\},$ 
and any Levi subgroup is given by $pL_sp^{-1}$ for $p\in P_s$. Their Lie algebras are given by
\begin{align*}
\mathfrak{p}_s &=\{Y\in\lieh^\C\ :\ \Ad(e^{ts})Y\text{ is
bounded as }t\to\infty\},\\
\mathfrak{l}_s &=\{Y\in\lieh^\C\ :\ \ad(Y)(s)=[Y,s]=0\}.
\end{align*}

We  consider the subspaces\newnot{fms}
\begin{align*}
&\liem_s=\{Y\in \liem^\C\ :\ \Ad(e^{ts})Y
\text{ is bounded as}\;\; t\to\infty\}\\
&\liem^0_{s}=\{Y\in \liem^\C\ :\ \Ad(e^{ts})Y=Y\;\;
\mbox{for every} \;\; t\}.
\end{align*}

One has that $\liem_s$ is invariant under the action of $P_{s}$ and $\liem^0_{s}$
is invariant under the action of $L_{s}$. 

\begin{remark}
 The subalgebra $\liem_s$ is the non-compact part of the parabolic subalgebra 
of $\lieg^\C$ defined by $s\in i\lieh$. Define 
$\widetilde{\liep}_s=\{ Y\in \lieg^\C \st \Ad(e^{ts})Y \textrm{ is bounded as }
t\to \infty \}.$  We have that $\liep_s = \widetilde{\liep}_s \cap \lieh^\C$ and 
$\liem_s = \widetilde{\liep}_s \cap \liem^\C$. Analogously, define 
$$\widetilde{\liel}_s=\{ Y\in \lieg^\C \st \Ad(e^{ts})Y = Y\;\; \mbox{for every} \;\; t\ \}.$$
Then $\liel_s = \widetilde{\liel}_s \cap \lieh^\C$ and $\liem^0_s = \widetilde{\liel}_s \cap \liem^\C$.
\end{remark}

An  element $s\in i\lieh$ defines a character $\chi_s$ of $\liep_s$ since
$\langle s,[\liep_s,\liep_s]\rangle=0$.  Conversely, by the isomorphism 
$\left( \liep_s/ [\liep_s,\liep_s]\right)^* \cong \liez_s^*$, 
where $\liez_s$ is the centre of the Levi subalgebra $\liel_s$, a 
character $\chi$ of $\liep_s$ is given by an element in $\liez_s^*$, 
which gives, via the invariant metric, an element 
$s_\chi\in \liez_s\subset i\lieh$. When $\liep_s\subset \liep_{s_\chi}$, 
we say that $\chi$ is an antidominant character of $\liep$. 
When  $\liep_s=\liep_{s_\chi}$ we say that 
$\chi$ is a strictly antidominant character. Note that for $s\in i\lieh$, 
$\chi_s$ is a strictly {\bf antidominant character} of $\liep_s$.

Let now $(E,\varphi)$ be a $G$-Higgs bundle over $X$, 
 and let $s\in i\lieh$. Let $P_s$ be defined as above.
For $\sigma\in \Gamma(E(H^\C/P_s))$ a reduction of the structure group 
of $E$ from $H^\C$ to $P_s$, we define the degree relative to $\sigma$ and
$s$, or equivalently to $\sigma$ and $\chi_s$, as follows. When  
a real multiple $\mu \chi_s$ of the character exponentiates to a 
character $\tilde{\chi}_s$ of $P_s$, 
we compute the degree as 
$$\deg(E)(\sigma,s)= \frac{1}{\mu} \deg(E_\sigma(\tilde{\chi}_s)).$$
This condition is not always satisfied, but one shows 
(\cite[Sec. 4.6]{GGM09}) 
that the antidominant character can be expressed as a linear combination of 
characters of the centre and fundamental weights, and there exists an integer 
multiple $m$ of the characters of the centre and the 
fundamental weights exponentiating to the group, so we can define the degree.

An alternative definition of the degree can be given in terms of the 
curvature of connections using Chern--Weil theory.
This definition is  more natural when considering gauge-theoretic 
equations as we do below. For this, define $H_s=H\cap L_s$ and
$\lieh_s=\lieh\cap\liel_s$. 
Then $H_s$ is a maximal compact subgroup of $L_s$, so the inclusions
$H_s\subset L_s$ is a homotopy equivalence. Since the inclusion
$L_s\subset P_s$ is also a homotopy equivalence, given a reduction
$\sigma$ of the structure group of $E$ to $P_s$ one can
further restrict the structure group of $E$ to $H_s$ in a unique
way up to homotopy. Denote by $E'_{\sigma}$ the resulting $H_s$
principal bundle.
Consider now a connection $A$
on $E'_{\sigma}$ and let
$F_A\in\Omega^2(X,E'_{\sigma}(\lieh_s)$ be  its
curvature. Then $\chi_s(F_A)$ is a $2$-form on $X$ with
values in $i\R$, and 
\begin{equation}\label{degree-chern-weil}
\deg(E)(\sigma,s):=\frac{i}{2\pi}\int_X \chi_s(F_A).
\end{equation}

We define the subalgebra $\lieh_{\ad}$\newnot{fhiota} as follows. 
Consider the decomposition $ \lieh = \liez + [\lieh, \lieh] $, where $\liez$ 
is the centre of $\lieh$, and the isotropy representation 
$\ad= \ad:\lieh\to \End(\liem)$. Let $\liez'=\ker(\ad_{|\liez})$  and 
take $\liez''$ such that $\liez=\liez'+\liez''$. Define the subalgebra 
$\lieh_{\ad} := \liez'' + [\lieh, \lieh]$. The subindex $\ad$ denotes that 
we have taken away the part of the centre $\liez$ acting trivially via 
the isotropy representation $\ad$.

\begin{remark}\label{remark:fz'-Hermitian-type}
  For groups of Hermitian type, $\liez'=0$ since an element both in 
$\liez$ and $\ker(\ad)$ belongs to the centre of $\lieg$, which is 
zero, as $\lieg$ is semisimple. Hence $\lieh_{\ad} = \lieh$.
\end{remark}

With $L_s$ $\liem_s$ and $\liem^0_s$ defined as above. We have the following.

\begin{definition}
\label{def:L-twisted-pairs-stability} 
Let  $\alpha\in i\liez\subset\liez^\C$.
 We say that a $G$-Higgs bundle $(E,\varphi)$ is:

$\alpha$-{\bf semistable} if for any $s\in i\lieh$ and any holomorphic reduction 
$\sigma\in\Gamma(E(H^\C/P_s))$ such that $\varphi\in
 H^0(X,E_{\sigma}(\liem_s)\otimes K)$, 
we have that $\deg(E)(\sigma,s)-\chi_s(\alpha)\geq 0.$

 $\alpha$-{\bf stable} if for any $s\in i\lieh_{\ad}$ and any holomorphic
  reduction $\sigma\in\Gamma(E(H^\C/P_s))$ such that $\varphi\in
  H^0(X,E_{\sigma}(\liem_s)\otimes K)$, we have that
  $\deg(E)(\sigma,s)-\chi_s(\alpha) > 0.$

$\alpha$-{\bf polystable} if it is $\alpha$-semistable and for
any $s\in i\lieh_{\ad}$ and any holomorphic reduction 
$\sigma\in\Gamma(E(H^\C/P_s))$ 
such that $\varphi\in  H^0(X,E_{\sigma}(\liem_s)\otimes K)$ and  
$\deg(E)(\sigma,s)-\chi_s(\alpha)=0$, there is a holomorphic reduction 
of the structure group $\sigma_L\in\Gamma(E_{\sigma}(P_s/L_s))$ to a Levi
subgroup $L_s$ such that  $\varphi\in H^0(X,E_{\sigma_L}(\liem_s^0)\otimes K)
\subset H^0(X,E_{\sigma}(\liem_s)\otimes K)$.
\end{definition}

\begin{remark}\label{remark:grupo-GL-Levi}
We may define a real group $G_{L_s} = (L_s\cap H) \exp(\liem^0_s\cap
\liem)$\newnot{GLS}  with maximal compact subgroup the compact real 
form $L_s\cap H$ of the complex group $L_s$ and $\liem^0_s\cap \liem$ 
as isotropy representation. Thus, an $\alpha$-polystable $G$-Higgs 
bundle reduces to an $\alpha$-polystable $G_{L_s}$-Higgs bundle 
since $\varphi$ belongs 
$H^0(X,E_{\sigma_L}(\liem_s^0)\otimes K)$. 
\end{remark}

\begin{remark}\label{twisted}
We can  replace $K$ in the definition of $G$-Higgs bundle by any holomorphic
line bundle $L$ on $X$. More precisely, an
\textbf{$L$-twisted $G$-Higgs bundle} $(E,\varphi)$ consists of a principal
$H^\C$-bundle $E$, and a holomorphic section $\varphi\in H^0(X,E(\liem^\C)\ot
L)$. We reserve the name $G$-Higgs bundle for the $K$-twisted case. 
The stability criteria  are as in Definition 3.5, replacing $K$ by $L$.
\end{remark}

\begin{remark}
For $G$ semisimple, the notion of $\al$-stability with $\al\neq 0$ only makes
sense for groups of Hermitian type, since $\al$ belongs to the centre of
$\lieh$, which is not zero if and only if the centre of a maximal compact
subgroup $H$ is non-discrete, i.e., if $G$ is of Hermitian type. For reductive
groups which are not of Hermitian type, $\al$-stability makes sense, but there
is only one value of $\al$ for which the  condition is not void. This value is
fixed by the topology of the principal bundle (see
\cite{garcia-prada-oliveira} for details).
\end{remark}

\begin{remark}
When $H^\C$ is connected, as mentioned above, every parabolic subgroup of
$H^\C$ is conjugate to one of the form $P_s$ for $s\in i\lieh$. In this 
situation, we can formulate the stability conditions in Definition 
\ref{def:L-twisted-pairs-stability} in terms of any parabolic subgroup 
$P\subset H^\C$, replacing $s$ by $s_\chi$, for any antidominant character 
$\chi$ of $\lieh^\C$.
\end{remark}

Two $G$-Higgs bundles $(E,\varphi)$ 
and $(E',\varphi')$ are isomorphic if there is an isomorphism $f:E\to E'$ such 
that $\varphi'=f^*\varphi$, where $f^*$ is the map $E(\liem^\C)\ot K \to E'(\liem^\C)\ot
K$  induced by $f$.

The {\bf moduli space of $\alpha$-polystable $G$-Higgs bundles} 
$\cM^\alpha(G)$  is defined as the set of isomorphism classes of 
$\alpha$-polystable  $G$-Higgs bundles on $X$. When $\alpha=0$ 
we simply denote  $\cM(G):= \cM^0(G)$.

\begin{remark}
Similarly, we can define the {\bf moduli space of $\alpha$-polystable 
$L$-twisted $G$-Higgs bundles} which will be denoted  by $\cM^\alpha_L(G)$.
\end{remark}

These moduli spaces  have the structure of a complex analytic
varieties,  as one can see by the standard slice method, which gives local 
models via the so-called Kuranishi map (see, e.g., \cite{kobayashi}).  
When $G$ is algebraic and under fairly general conditions, 
the moduli spaces $\cM^\al(G)$ can be constructed by geometric invariant 
theory and hence are complex algebraic varieties. The  work of 
Schmitt \cite{schmitt08} deals with the construction 
of the moduli space of $L$-twisted $G$-Higgs bundles for $G$ a real reductive 
Lie group. This construction generalizes the constructions of the moduli space 
of $G$-Higgs bundles done by Ramanathan \cite{ramanathan96} when $G$ is
compact ($G$-Higgs bundles in this case are simply $G^\C$-principal bundles),  
and by Simpson \cite{simpson94,simpson95} when $G$ 
is a complex reductive algebraic (see also \cite{nitsure} for $G=\GL(n,\C)$).

\begin{remark}\label{gln}
If  $G=\GL(n,\C)$, we recover the
original notions of Higgs bundle and stability
introduced by Hitchin \cite{hitchin1987}.
A  $G=\GL(n,\C)$-Higgs bundle is given by
a holomorphic vector bundle $V:=E(\C^n)$ --- associated to a principal
$\GL(n,\C)$-bundle $E$ via the standard representation ---
and a homomorphism
$$
\Phi: V\lra V\otimes K.
$$
The Higgs bundle $(V,\Phi)$ is   stable if
\begin{equation}\label{higgs-stability}
\frac{\deg V'}{\rank V'}< \frac{\deg V}{\rank V}
\end{equation}
for every proper subbundle $V'\subset V$ such that
$\Phi(V')\subset V'\otimes K$.
The Higgs bundle $(V,\Phi)$ is
polystable if $(V,\Phi)=\oplus_i (V_i,\Phi_i)$ where
$(V_i,\Phi_i)$ is a  stable Higgs bundles  and
${\deg V_i}/{\rank V_i}= {\deg V}/{\rank V}$.

If we consider $G=\SL(n,\C)$ we have to add  that 
$\det(V)\cong \cO$ (hence $\deg V=0$) and $\Tr(\Phi)=0$.

\end{remark}

The notion of stability emerges from the study of the
Hitchin equations.  The equivalence between the existence of solutions to 
these equations and the $\al$-polystability of Higgs bundles is known as 
the {\bf Hitchin-Kobayashi correspondence}, which we state below.

\begin{theorem}\label{theo:hk-twisted-pairs}
Let $(E,\varphi)$ be a $G$-Higgs bundle over a Riemann surface $X$ with volume
form $\omega$. Then $(E,\varphi)$ is $\alpha$-polystable if and only if 
there exists a reduction $h$ of the structure group of $E$ from $\H^\C$ to
$H$, that is a smooth section of $E(H^\C/H)$,  such that
\begin{equation}\label{eq:Hitchin-Kobayashi-h}
F_h - [\varphi,\tau_h(\varphi)]=-i\al \omega\\
\end{equation}
where $\tau_h:\Omega^{1,0}(E(\liem^\C))\to \Omega^{0,1}(E(\liem^\C))$ is the
combination of  the anti-holomorphic involution in $E(\liem^\C)$ defined 
by the compact real form at each point determined by $h$  and the 
conjugation of $1$-forms, and $F_h$ is the 
curvature of the unique $H$-connection compatible with the holomorphic 
structure of $E$ (the Chern connection). 
\end{theorem}

This theorem was proved by Hitchin in the case of $\SL(2,\C)$, by Simpson when
$G$ is complex, and in \cite{BGM03,GGM09} for a general 
reductive real Lie group $G$

\begin{remark}\label{theo:-twisted-hk-twisted-pairs}
There is a theorem similar  to Theorem \ref{theo:hk-twisted-pairs} for
$L$-twisted $G$-Higgs bundles (see Remark \ref{twisted} and
\cite{GGM09}) 
for an arbitrary line bundle $L$. If $(E,\varphi)$ is such a pair,
one fixes a Hermitian metric $h_L$ on $L$, and looks for
a reduction of structure group $h$  of $E$ from $H^\C$ to $H$ satisfying

\begin{equation}\label{eq:twisted-Hitchin-Kobayashi-h}
F_h - [\varphi,\tau_h(\varphi)]\omega=-i\al \omega,
\end{equation} 
where now 
$\tau_h:\Omega^0(E(\liem^\C)\otimes L)\to \Omega^0(E(\liem^\C)\otimes L)$ is 
the combination of the anti-holomorphic involution in $E(\liem^\C)$ defined
by the compact real form at each point determined by $h$  and the 
metric $h_L$.
\end{remark}

\subsection{Higgs bundles and representations}\label{higgs-reps}

Fix a base point $x\, \in\, X$.
A \textbf{representation} of $\pi_1(X,x)$ in
$G$ is  a homomorphism $\pi_1(X,x) \,\longrightarrow\, G$.
After fixing a presentation of $\pi_1(X,x)$, the set of all such homomorphisms,
$\Hom(\pi_1(X,x),\, G)$, can be identified with the subset
of $G^{2g}$ consisting of $2g$-tuples
$(A_{1},B_{1}, \cdots, A_{g},B_{g})$ satisfying the algebraic equation
$\prod_{i=1}^{g}[A_{i},B_{i}] \,=\, 1$. This shows that
$\Hom(\pi_1(X,x),\, G)$ is an algebraic variety.

The group $G$ acts on $\Hom(\pi_1(X,x),G)$ by conjugation:
$$
(g \cdot \rho)(\gamma) \,=\, g \rho(\gamma) g^{-1}\, ,
$$
where $g \,\in\, G$, $\rho \,\in\, \Hom(\pi_1(X,x),G)$ and
$\gamma\,\in \,\pi_1(X,x)$. If we restrict the action to the subspace
$\Hom^+(\pi_1(X, x),\,G)$ consisting of reductive representations,
the orbit space is Hausdorff.  We recall that a \textbf{reductive representation}
is one whose composition with the adjoint representation in $\mathfrak g$
decomposes as a direct sum of irreducible representations.
This is equivalent to the condition that the Zariski closure of the
image of $\pi_1(X,x)$ in $G$ is a reductive group. Define the
{\bf moduli space of representations} of $\pi_1(X,x)$ in $G$ to be the orbit space
$$
\mathcal{R}(G) = \Hom^{+}(\pi_1(X,x),G) /G.
$$
This is a real algebraic variety.
For another point $x'\, \in\, X$, the fundamental groups
$\pi_1(X,x)$ and $\pi_1(X,x')$ are identified by an isomorphism unique up to
an inner automorphism. Consequently, $\mathcal{R}(G)$ is independent of the choice of
the base point $x$.

Given a representation $\rho\colon\pi_{1}(X,x) \,\longrightarrow\,
G$, there is an associated flat principal $G$-bundle on
$X$, defined as
$$
  E_{\rho} \,=\, \widetilde{X}\times_{\rho}G\, ,
$$
where $\widetilde{X} \,\longrightarrow\, X$ is the universal cover
and $\pi_{1}(X, x)$ acts
on $G$ via $\rho$.
This gives in fact an identification between the set of equivalence classes
of representations $\Hom(\pi_1(X),G) / G$ and the set of equivalence classes
of flat principal $G$-bundles, which in turn is parametrized by
the (non-abelian) cohomology set $H^1(X,\, G)$. We have the following.

\begin{theorem}\label{na-Hodge}
Let $G$ be a semisimple real Lie group. Then there is a homeomorphism
$\mathcal{R}(G) \,\cong\, \cM(G)$. 
\end{theorem}

The proof of Theorem \ref{na-Hodge} is the combination of Theorem
\ref{theo:hk-twisted-pairs} and the following  theorem of Corlette \cite{corlette},
also proved by Donaldson \cite{donaldson} when $G=\SL(2,\C)$.

\begin{theorem}\label{corlette}
Let $\rho$ be a representation of $\pi_1(X)$ in $G$ with corresponding
flat $G$-bundle $E_\rho$. Let $E_\rho(G/H)$ be the associated $G/H$-bundle.
Then the existence of a harmonic section of $E_\rho(G/H)$ is equivalent to
the reductiveness of $\rho$. 
\end{theorem}

\begin{remark}
Theorem \ref{na-Hodge} can be extended to reductive groups replacing 
$\pi_1(X)$ by its universal central extension.
\end{remark}

We can  assign  a topological invariant  to a representation
$\rho$ given by the  characteristic class $c(\rho):=c(E_{\rho})\in \pi_1(G)$
corresponding
to $E_{\rho}$. To define this, let $\widetilde G$ be the universal covering group
of $G$. We have an exact  sequence
$$
1 \lra \pi_1(G)\lra \widetilde G \lra G \lra 1
$$
which gives rise to the (pointed sets) cohomology sequence
\begin{equation}\label{characteristic}
H^1(X, {\widetilde G}) \lra H^1(X, {G})\stackrel{c}{\lra}   H^2(X, \pi_1(G)).
\end{equation}

Since $\pi_1(G)$ is abelian, we have
$$
 H^2(X, \pi_1(G))\cong  \pi_1(G),
$$
and $c(E_\rho)$ is defined as the image of $E$ under the last map in
(\ref{characteristic}). Thus the class $c(E_\rho)$ measures  the
obstruction to lifting $E_\rho$ to a flat $\widetilde G$-bundle, and hence to
lifting $\rho$ to a representation of $\pi_1(X)$ in $\widetilde G$.  
For a  fixed
 $c\in \pi_1(G)$, 
we can consider the subvariety of the  moduli space of representations
 $\mathcal{R}_c(G)$ defined by 
$$
\mathcal{R}_c(G):=\{\rho \in \mathcal{R}(G)\; : \; c(\rho)=c\}.
$$

Similarly, we can assign a topological invariant to a $G$-Higgs bundle 
$(E,\varphi)$.
Assuming that $G$ is connected (and hence $H$ is connected),
topologically, $H^\C$-bundles $E$  on $X$ are classified by
a characteristic class $c=c(E)\in \pi_1(H^\C)=\pi_1(H)=\pi_1(G)$. Again, this 
comes from considering  the exact sequence  
$$
1 \lra \pi_1(H^\C)\lra \widetilde{H^\C} \lra H^\C \lra 1
$$
and the corresponding  (pointed sets) cohomology sequence
$$
H^1(X, \underline{{\widetilde H^\C}}) \lra  H^1(X, \underline{H^\C})
\stackrel{c}{\lra}   H^2(X, \pi_1(H^\C)).
$$
For a fixed  such  class  $c$, the subvariety  $\mathcal{M}_c(G)\subset \cM(G)$
is defined as  the set of isomorphism classes of polystable
$G$-Higgs bundles $(E,\varphi)$  such that $c(E)=c$.

If $G$ is semisimple, of course, the homeomorphism in Theorem \ref{na-Hodge} restricts to give a homeomorphism
$$
\mathcal{R}_c(G) \,\cong\, \cM_c(G). 
$$

\begin{proposition}\label{components}
Let $G$ be a connected semisimple real Lie group.  Then

(1) There is a homomorphism
$$
\pi_0(\mathcal{R}(G)) \,\cong\, \pi_0(\cM(G))\lra \pi_1(G). 
$$

(2) If $G$ is compact or complex, the homomorphism in (1) is an isomorphism.
 \end{proposition}

The proof of (2) in Proposition \ref{components} is due J. Li \cite{li}. An 
alternative proof using Higgs bundles, and including the case of non-connected 
reductive Lie groups is given in \cite{garcia-prada-oliveira}.

The homomorphism in (1) of Proposition \ref{components} for a real form of a complex 
semisimple Lie group may be neither surjective nor injective, as we will see in Section 
\ref{sl2r}.

\section{A basic example: $G=\SL(2,\R)$}\label{sl2r}

We consider $G=\SL(2,\R)$. Let $X$ be a compact Riemann surface of genus $g\geq 2$.
Given a representation $\rho:\pi_1(X)\to \SL(2,\R)$, the invariant associated to $\rho$, 
defined in Section \ref{higgs-reps} is  an integer $d \in\pi_1(\SL(2,\R)\cong \Z$. 
This is the Euler class $d\in\Z$ of the corresponding flat $\SL(2,\R)$-bundle $E_\rho$.
We can then define the subvarieties
$$
\cR_d:=\{\rho\in \cR(\SL(2,\R))\;:\; \mbox{with Euler class} \;\;d\}.
$$

By Milnor~\cite{milnor}, $\cR_d$ is empty unless 
\begin{equation}\label{milnor}
|d|\leq g-1.
\end{equation}

A maximal compact subgroup of $\SL(2,\R)$ can be identified with $\U(1)$, and 
its complexification $\C^*$ can be embedded in  $\SL(2,\C)$ as the matrices of the form 
$$
\begin{pmatrix}
  a & 0 \\
  0& a^{-1}
\end{pmatrix},
$$
with $a\in \C^*$. On the other hand, in this case $\liem^\C$ coincides with the matrices
$$
\begin{pmatrix}
  0 & b \\
  c  & 0
\end{pmatrix},
$$
with $b,c\in \C$.

An $\SL(2,\R)$-Higgs bundle is then simply a triple $(L,\beta,\gamma)$,
where $L$ is a holomorphic line bundle over $X$ and $\beta\in H^0(X,L^2K)$ and $\gamma\in H^0(X,L^{-2}K)$. 
It can also be viewed as the 
$\SL(2,\C)$-Higgs bundle $(V,\Phi)$ obtained by extension of structure group given by
$V=L\oplus L^{-1}$, and 
$$
\Phi =
\begin{pmatrix}
  0 & \beta \\
  \gamma  & 0
\end{pmatrix}.
$$

The topological invariant attached to the Higgs bundle  $(L,\beta,\gamma)$ is then $\deg L$, the degree of $L$ (its first Chern class).
As for representations, we can  define
the subvariety  $\cM_d\subset \cM(\SL(2,\R))$ as the  $\SL(2,\R)$-Higgs bundles  with fixed $\deg
L=d$.  Here the natural gauge transformations are
$\C^*$-transformations, that is, those of $L$. Allowing
$\SL(2,\C)$-transformation naturally  identifies $\cM_d$ and
$\cM_{-d}$  inside $\cM(\SL(2,\C))$.
Now, the semistability condition gives a constraint
on the possible degrees that $L$ may have, namely,
we must have
\begin{equation}\label{milnor-higgs}
\abs{d}\leq g-1.
\end{equation}
This can be seen very easily (\cite{hitchin1987}). Assume that $d\geq 0$ (similar
argument for $d\leq 0$). Suppose that $d>g-1$. Then $\gamma=0$ and
$L$ is a $\Phi$-invariant line subbundle of $V$. By the semistability
of  $(V,\Phi)$ we must have that $d\leq 0$ which gives a contradiction.
Considering the homeomorphism $\cR_d\cong\cM_d$, this gives a Higgs bundle proof of the 
Milnor inequality \cite{milnor}.

The  moduli space of representations of $\pi_1(X)$ in $\SL(2,\R)$
was studied by Goldman \cite{goldman}, who showed
that for $d$ satisfying $d=g-1$ (same for $d=1-g$) there are $2^{2g}$ isomorphic
connected components that can be identified with Teichm\"uller
space, and showed that for $d$ such that  $\abs{d}<g-1$ there is
only one connected component.  
This was also proved by Hitchin \cite{hitchin1987}, who also gave a very explicit description of each
component. Namely, for $0<|d|<g-1$, $\cM_d$ is the total space of a holomorphic  complex vector bundle
of rank $g+2|d|-1$ over a $2^{2g}$-fold covering of the
$2g-2-2|d|$-symmetric power of $X$
$$
\Sym^{2g-2-2|d|}X.
$$
To see this, assume  that $d>0$ (the case $d< 0$ is similar), then the pair
$(L,\gamma)$ defines an element in $\Sym^{2g-2-2d}(X)$, given by the
zeros of $\gamma$. But if $L_0$ is a line bundle of degree $0$ of
order two, that is, $L_0^2=\cO$, then the element $(LL_0,\gamma)$
defines also the same divisor in $\Sym^{2g-2-2d}(X)$. Hence the set 
of pairs $(L,\gamma)$ gives a point in the $2^{2g}$-fold covering of
$\Sym^{2g-2-2d}(X)$. The section $\beta$ now gives the fibre of the
vector bundle. Notice  Riemann--Roch implies that  $H^1(X,L^2K)=0$.
The case of $d=0$ is a bit more involved (see \cite{hitchin1987} and \cite{garcia-prada-ramanan})

From this we conclude the following.

\begin{proposition}

(1)  $\dim \cM_d=3g-3$,

(2) $\cM_d$ is connected if $|d|< g-1$,

(3) $\cM_d$ has $2^{2g}$ connected components if $|d|=g-1$, each
isomorphic  to $\C^{3g-3}$ (the fibre of a rank $3g-3$ vector bundle
over a $2^{2g}$-fold covering of a point!). 
\end{proposition}

(3) is clear since if $\deg L=g-1$, the line bundle $L^{-2}\otimes K$ is of zero 
degree and hence has a section (unique up to multiplicaton by a scalar) if
and only if $L^{-2}\otimes K\cong \cO$, i.e. if $L$ is a square root of $K$.
For each of the $2^{2g}$ choices of square root $L=K^{1/2}$, one has a connected 
component  which is parametrized by $\beta\in H^0(X,K^2)$. 
Each of these components is 
diffeomorphic to Teichm\"uller space.

\section{The Hitchin map, split real forms and Hitchin components}
For details on this section, see \cite{GPR18}.

\subsection{Maximal split subgroup and Chevalley map}\label{kostant-rallis}

Let $\lieg$ be a reductive real Lie algebra with a Cartan involution 
$\theta$ decomposing
$\lieg$ as $\lieg={\lieh}\oplus\liem.$
Given a maximal subalgebra 
$\liea\subset\liem$ it follows that it must be abelian, and one can easily 
prove that its elements are semisimple and diagonalizable over the real numbers (cf. \cite[Chap.VI]{knapp}, 
note that
Knapp proves it for semisimple Lie algebras, but for reductive Lie algebras it suffices to use invariance of the
centre
and the semisimple part of $[\lieg,\lieg]$) under the Cartan involution.
Any such subalgebra is called a \textbf{maximal anisotropic} Cartan subalgebra of $\lieg$. By extension,
its complexification $\liea^\C$ is called a \textbf{maximal anisotropic} Cartan subalgebra of $\lieg^\C$ (with respect
to $\lieg$). A
maximal anisotropic Cartan subalgebra $\liea$ can be completed to a $\theta$-equivariant \textbf{Cartan subalgebra} of $\lieg$,
namely, a subalgebra whose complexification is a Cartan subalgebra of $\lieg^\C$. Indeed, define
$$
\lied=\liet\oplus\liea 
$$
where $\liet\subset\liec_{\lieh}(\liea):=\{x\in\lieh\ :\ [x,\liea]=0\}$ is a maximal abelian subalgebra 
(\cite{knapp}, Proposition 6.47). Cartan subalgebras of this kind
(and their complexifications) are called \textbf{maximally split}.

The dimension of maximal anisotropic Cartan subalgebras of a real reductive Lie algebra $\lieg$ is called the \textbf{the
real (or split) rank of }
$\lieg$. This number measures the degree of compactness of real forms: indeed,  a real form 
is compact (that is, its adjoint group is compact)
 if and only if $\mathrm{rk}_\R(\lieg)=0$. On the other hand, a real form is defined to
be \textbf{split} 
if $\mathrm{rk}_\R(\lieg)=\mathrm{rk}{\lieg^\C}$. Note that the split rank depends on the involution $\theta$ associated with the real form, when $\lieg$ is not semisimple.

 In \cite{KR71}, Kostant and Rallis give a procedure to construct a $\theta$-invariant subalgebra 
$\tlieg\subset\lieg$ such $\tlieg\subset(\tlieg)^\C$ is a split real form, whose Cartan subalgebra is 
${\liea}$ and such that $\liez(\tlieg)=\liez(\lieg)\cap\liem$. Their construction relies on the following notion. A \textbf{three dimensional subalgebra} (TDS) 
$\lies^\C\subset{{\lieg^\C}}$ 
is the image of 
an injective morphism $\liesl(2,\C) \to{\lieg^\C}$.
 A TDS is called \textbf{normal} if  $\dim \lies^\C\cap{\lieh^\C}=1$ and 
$\dim \lies^\C\cap{\liem^\C}=2$.
 It is called \textbf{principal} if it is generated by elements $\{e,f,x\}$, 
where $e$ and $f$ are nilpotent regular elements in $\liem^\C$,  
and $\ {x}\in\lieh^\C$ is semisimple. 
A set of generators satisfying such relations is called a
\textbf{normal triple}. 

A subalgebra $\tlieg\subset\lieg$ generated by $\liea$ and $\lies^\C\cap \lieg$, 
where $\lies^\C$ is a 
principal normal TDS invariant by the involution defining $\lieg$ inside of 
$\lieg^\C$ is called a \textbf{maximal split subalgebra}.

\begin{remark}\label{split form complex}
Let ${\lieg^\C}$ be a complex reductive Lie algebra, and let $\left({\lieg^\C}\right)_\R$ be its underlying real reductive
algebra. Then,
the maximal split subalgebra of $\left({\lieg^\C}\right)_\R$ is isomorphic to the split real form ${\lieg}_{\Split}$ of ${\lieg^\C}$. 
It is clearly split within its complexification and it is maximal within $\left({\lieg^\C}\right)_\R$ with this property,
which can be easily checked by identifying $\left({\lieg^\C}\right)_\R\cong\lieg_{\Split}\oplus i\lieg_{\Split}$.
\end{remark}

We say that a real Lie group is split if its Lie algebra is split.
Let  $(G,H,\theta, B)$ be a real reductive Lie group.
The tuple $(\widehat{G},\widehat{H},\widehat{\theta},\widehat{B})$,
where $\widehat{G}\subset G$ is the analytic subgroup correspoding to 
the maximal split subalgebra $\widehat{\lieg}\subset \lieg$, $\widehat{H}:=
\exp(i\widehat{\lieh})\leq H$, and $\widehat{\theta}$ 
and $\widehat{B}$ are obtained by restriction,  is a  reductive Lie 
group to which we refer as  the {\bf maximal split subgroup} of   
$(G,H,\theta, B)$.

The restriction to $\liea$ of the adjoint representation of $\lieg$ yields 
a decomposition of $\lieg$ into  $\liea$-eigenspaces
$$\lieg=\bigoplus_{\lambda\in\Lambda(\liea)}{\lieg}_{\lambda},$$
where $\Lambda(\liea)\subset\liea^{*}$ is called the set of 
\textbf{restricted roots} of $\lieg$ with respect to $\liea$.
The set $\Lambda({\liea})$ forms a root system 
(see \cite[Chap. II, Sec. 5]{knapp}), which may not
be reduced (that is, there may be roots whose double is also a root).
The name restricted roots is due to the following fact:  extending restricted 
roots by $\C$-linearity, we obtain 
$\Lambda({\liea^\C})\subset \left(\liea^\C\right)^*$, also called 
restricted roots.

We define the \textbf{restricted Weyl group} of   ${\lieg^\C}$ 
associated to  ${\liea^\C}$, $W({\liea^\C})$, to be
 the group of automorphisms of ${\liea^\C}$ 
generated by  reflections on the hyperplanes defined
 by the restricted roots $\lambda\in\Lambda(\liea^\C)$.
\begin{proposition}\label{prop exponents}
Let $(G, H,\theta,B)$ be a real  reductive Lie group and 
$(\widehat{G},\widehat{H},\widehat{\theta},\widehat{B})<(G, H,\theta,B)$ be 
its  maximal  split subgroup.  Let $a=\dim\liea^\C$. Then:

\begin{enumerate}
\item  Restriction induces 
an isomorphism
 $$
 \C[{\liem^\C}]^{{H^\C}}\cong\C[{\liea^\C}]^{W({\liea^\C})}
\cong\C[\widehat{{\liem}}^\C]^{\widehat{H}^\C}.
 $$
 
\item $\C[\liea^\C]^{{W(\liea^\C)}}$ is generated by homogeneous polynomials 
of degrees $d_1,\dots, d_a$, canonically determined
by $(G,\theta)$.  

\item The degrees of the homogeneous polynomials for $(G, H,\theta,B)$ and 
$(\widehat{G},\widehat{H},\widehat{\theta},\widehat{B})$ coincide
 \end{enumerate}
\end{proposition}

We thus have an algebraic morphism, called the {\bf Chevalley map}
\begin{equation}\label{Chevalley map}
\chi: {\liem^\C}\projects {\liem^\C}\sslash{H^\C}\cong {\liea^\C}/W({\liea^\C})
\end{equation}
where the double quotient sign $\sslash$ stands for the affine GIT quotient.

\subsection{The Hitchin map}
 
Let $(G,H,\theta,B)$ be a reductive real Lie group. Consider the Chevalley 
morphism  (\ref{Chevalley map}):
$$
\chi:\liem^\C\to{\liea^\C}/W({\liea^\C}).
$$
This map is  $\C^\times$-equivariant. In particular, it induces a morphism
$$
h:\liem^\C\otimes K\to{\liea^\C}\otimes K/W({\liea^\C}).
$$
The map $\chi$ is also $H^\C$-equivariant, thus defining a morphism
\begin{equation}\label{defi stacky Hitchin map}
h: \cM(G)\to B(G):=H^0(X,{\liea^\C}\otimes K/W({\liea^\C})).
\end{equation}

The map $h$ is called the \textbf{Hitchin map}, and the space $B(G)$ is 
called the \textbf{Hitchin base}.

In more concrete terms Let $p$ be an $H^\C$-invariant homogeneous polynomial of degree $d$ on
$\liem^\C$. Then, if $(E,\varphi)$ is a  $G$-Higgs
bundle, the evaluation of $p$ on $\varphi$ gives an element of
$H^0(X,K^{d})$, so $p$ induces a map
$$
p:\mathcal{M}(G)\longrightarrow H^0(X,K^d).
$$
Let $p_1,\ldots,p_r$  with $r=\rk(G)$ be  a basis of the ring  of
$H^\C$-invariant homogeneous polynomials on $\liem^\C$. This basis defines a
map
$$
\mathcal{M}(G)\rightarrow\bigoplus_{i=1}^rH^0(X,K^{d_i})
$$
where, for $i=1,2,\ldots,r$, $d_i$ is the degree of $p_i$. These degrees are independent 
of the bases and they only depend on $G$. 
This is the Hitchin map  and hence the Hitchin base is 
$$
B(G)=\bigoplus_{i=1}^rH^0(X,K^{d_i}).
$$

If $G$ is complex and  $E$ is a stable $G$-bundle, one has 
(see \cite{hitchin:duke})
that
$$
\dim B(G)=\dim H^0(X,E(\lieg)\otimes K)
$$ 
and hence the dimension of $B(G)$ coincides with the dimension of the 
moduli space of $G$-bundles, which is half the dimension of $\cM(G)$. 
In \cite{hitchin:duke} Hitchin shows that $\cM(G)$ is a completely 
algebraically integrable system and that, for a classical group,
the generic fibre is either a Jacobian or a Prym variety of a curve covering 
$X$, which is called  the {\bf spectral curve}. The 
description of the fibres for a general complex reductive Lie group is
given by  Donagi and  Gaitsgory \cite{DG01} in terms of a {\bf cameral cover} 
of the curve. The description of the fibration for general real forms is not 
yet fully understood (for some partial  results see \cite{Peo13,Sch13,HS14,GP19}).

\subsection{Hitchin components for $\SL(n,\mathbb{R})$-Higgs bundles}
\label{split-slnr}

We will illustrate first how Hitchin constructed in \cite{hitchin92} special  
 components of the moduli space of $\SL(n,\R)$-Higgs bundles as the image of 
certain sections of the Hitchin map.
To explain  this, recall that Cartan decomposition of $\mathfrak{sl}(n,\mathbb{R})$ is given by
$$\mathfrak{sl}(n,\mathbb{R})=\mathfrak{so}(n)+\mathfrak{m},$$
where $\mathfrak{m}=\{\text{symmetric real matrices of trace
}0\}$. A $\SL(n,\mathbb{R})$- Higgs bundle is thus  a pair
$(E,\varphi)$,  where
$E$ is a principal holomorphic $\SO(n,\mathbb{C})$-bundle over $X$
and the Higgs field is a holomorphic section 
$$\varphi\in
H^0(E(\mathfrak{m}^\mathbb{C})\otimes K).
$$

Using the standard representations of $\SO(n,\mathbb{C})$ in
$\mathbb{C}^n$ one  can associate to $E$ a holomorphic vector bundle
$V$ of rank $n$ with $\det V=\mathcal{O}$ together with a
nondegenerate quadratic form $Q\in H^0(S^2 V^*)$.

A $\SL(n,\mathbb{R})$-Higgs bundle is then in correspondence with
a triple 
$$
(V,Q,\varphi),
$$
where the Higgs field is a symmetric and
traceless endomorphism $\varphi:V\rightarrow V\otimes K$.

The simplest case is to consider the complex Lie group
$\SL(2,\mathbb{C})$ and its split real form $\SL(2,\mathbb{R})$.
The Lie algebra $\mathfrak{sl}(2,\mathbb{C})$ has rank $1$ and the
algebra of invariant polynomials on it is generated by $p_2$ of
degree $2$ obtained from the characteristic polynomial
$$\det(x1-A)=x^2+p_2(A)$$ of a trace-free matrix. We are going to
define a section of the Hitchin map
\begin{eqnarray*}
  p:\mathcal{M} &\rightarrow& H^0(K^2) \\
  (E,\varphi) &\mapsto& p_2(\varphi)
\end{eqnarray*}
where here $\mathcal{M}$ denotes the moduli space of polystable
$\SL(2,\mathbb{C})$-Higgs bundles. This section will give an
isomorphism between the vector space $H^0(K^2)$ and a connected
component of the moduli space
$\mathcal{M}(\SL(2,\mathbb{R}))\subset\mathcal{M}$ of polystable
$\SL(2,\mathbb{R})$-Higgs bundles. To construct the section, we
consider the elements $$\big\langle x=\left(%
\begin{array}{cc}
  1 & 0 \\
  0 & -1 \\
\end{array}%
\right),e=\left(%
\begin{array}{cc}
  0 & 1 \\
  0 & 0 \\
\end{array}%
\right),\tilde{e}=\left(%
\begin{array}{cc}
  0 & 0 \\
  1 & 0 \\
\end{array}%
\right)\big\rangle\cong\mathfrak{sl}(2,\mathbb{C})$$ that satisfy
$$[x,e]=2e\text{, }[x,\tilde{e}]=-2\tilde{e}\text{ and
}[e,\tilde{e}]=x,$$ where $x$ is an element of the Cartan
subalgebra (a semisimple element) and $e$, $\tilde{e}$ are
nilpotent. The pair $(K^{1/2}\oplus K^{-1/2},\varphi=\tilde{e}-\alpha e=\left(%
\begin{array}{cc}
  0 & -\alpha \\
  1 & 0 \\
\end{array}%
\right))$, where $\alpha\in H^0(K^2)$, is a
$\SL(2,\mathbb{R})$-Higgs bundle. In the vector bundle
$K^{1/2}\oplus K^{-1/2}$ we
have the orthogonal structure $Q=\left(%
\begin{array}{cc}
  0 & 1 \\
  1 & 0 \\
\end{array}%
\right)$ and the Higgs field is symmetric with respect to this
orthogonal form. The section is finally defined by $$s(\alpha)=
\left(K^{1/2}\oplus K^{-1/2},\varphi=\left(%
\begin{array}{cc}
  0 & -\alpha \\
  1 & 0 \\
\end{array}%
\right)\right).$$ That is, the pairs $$\left\{(K^{1/2}\oplus K^{-1/2},\varphi=\left(%
\begin{array}{cc}
  0 & -\alpha \\
  1 & 0 \\
\end{array}%
\right))\right\}_{\alpha\in H^0(K^2)}$$ form a connected component
of $\mathcal{M}(\SL(2,\mathbb{R}))$ of dimension $6g-6$, and there are
$2^{2g}$ connected components isomorphic to this one --- the number
of possible choices of the square $K^{1/2}$. These are precisely the 
components with maximal integer invariant given by the Milnor inequality,
described in Section \ref{sl2r}.

Now consider the general case $\SL(n,\mathbb{R})$ which is the
split real form of $\SL(n,\mathbb{C})$. The Lie algebra
$\mathfrak{sl}(n,\mathbb{C})$ has rank $n-1$ and a basis for the
invariant polynomials on $\mathfrak{sl}(n,\mathbb{C})$ is provided
by the coefficients of the characteristic polynomial of a
trace-free matrix,
$$\det(x-A)=x^n+p_1(A)x^{n-2}+\ldots+p_{n-1}(A),$$ where
$\deg(p_i)=i+1$. Ehe Hitchin map
$$h:\mathcal{M}(\SL(n,\mathbb{C}))\rightarrow\displaystyle{\bigoplus_{i=1}^{n-1}
H^0(K^{i+1})}$$  is defined by
$$h(E,\varphi)=(p_1(\varphi),\ldots,p_{n-1}(\varphi)),$$ where
$\mathcal{M}(\SL(n,\mathbb{C}))$ is the moduli space of polystable
$\SL(n,\mathbb{C})$-Higgs bundles.

We are going to define a section of this map that will give an
isomorphism between the vector space
$\displaystyle{\bigoplus_{i=1}^{n-1} H^0(K^{i+1})}$ and a
connected component of the moduli space
$\mathcal{M}(\SL(n,\mathbb{R}))\subset\mathcal{M}(\SL(n,\mathbb{C}))$
of polystable $\SL(n,\mathbb{R})$-Higgs bundles.

A nilpotent element $e\in\mathfrak{sl}(n,\mathbb{C})$ is called
regular if its centralizer is $(n-1)$-dimensional. In
$\mathfrak{sl}(n,\mathbb{C})$ a regular nilpotent element is
conjugate to an element
$$e=\sum_{\alpha\in\Pi}X_\alpha,$$ where $\Pi=\{\alpha_i=e_i-e_{i+1},1\leq i\leq n-1\}$ and $X_{\alpha_i}=E_{i,i+1}$
is a root vector for $\alpha_i$, that is, $$e=\small{\left(%
\begin{array}{cccccc}
  0 & 1 & 0 & \cdots & \cdots & 0 \\
  \vdots & 0 & 1 & 0 & \cdots & 0 \\
  \vdots & \vdots &  & \ddots &  & \vdots \\
  \vdots & \vdots &  &  & \ddots & 0 \\
  \vdots & \vdots &  &  &  & 1\\
  0 & 0 & \cdots & \cdots & \cdots & 0\\
\end{array}%
\right).}$$ Any nilpotent element can be embedded in a
$3$-dimensional simple subalgebra $\langle
x,e,\tilde{e}\rangle\cong\mathfrak{sl}(2,\mathbb{C})$, where $x$
is semisimple, $e$ and $\tilde{e}$ are nilpotent, and they satisfy
$$\xymatrix{[x,e]=2e;&[x,\tilde{e}]=-2\tilde{e};&[e,\tilde{e}]=x.}$$
The adjoint action $$\langle
x,e,\tilde{e}\rangle\cong\mathfrak{sl}(2,\mathbb{C})\rightarrow\End(\mathfrak{sl}(n,\mathbb{C}))$$
of this subalgebra breaks up the Lie algebra
$\mathfrak{sl}(n,\mathbb{C})$ as a direct sum of irreducible
representations
$$\mathfrak{sl}(n,\mathbb{C})=\bigoplus_{i=1}^{n-1} V_i,$$ with $\dim(V_i)=2i+1$. That is, each $V_i$ is the
irreducible representation $S^{2i}\mathbb{C}^2$, where
$\mathbb{C}^2$ is the standard representation of
$\mathfrak{sl}(2,\mathbb{C})$, and the eigenvalues of $\ad x$ on
$V_i$ are $-2i,-2i+2,...,2i-2,2i$.

The highest weight vector of $V_i$, defined as a vector $e_i\in
V_i$ that is an  eigenvector for the action of $x$ and is in the
kernel of $\ad(e)$, has eigenvalue $2i$ for $\ad x$. We take
$V_1=\langle x,e,\tilde{e}\rangle$ and $e=e_1$.

Given
$(\alpha_1,...,\alpha_{n-1})\in\displaystyle{\bigoplus_{i=1}^{n-1}
H^0(K^{i+1})}$, we define the Higgs field in the Hitchin component
by
$$\varphi=\tilde{e}_1+\alpha_1 e_1+...+\alpha_{n-1} e_{n-1},$$ and
the vector bundle is given by $$S^{n-1}(K^{-1/2}\oplus
K^{1/2})=K^{-(n-1)/2}\oplus K^{-(n-3)/2}\oplus\cdots\oplus
K^{(n-3)/2}\oplus K^{(n-1)/2}.$$ The field $\varphi$ is given in
the following order of the basis, $K^{(n-1)/2}\oplus
K^{(n-3)/2}\oplus\cdots\oplus K^{-(n-3)/2}\oplus K^{-(n-1)/2}$.

Hitchin \cite{hitchin92} 
generalized this construction to the split real form of any 
complex semisimple Lie group, proving the following. 

\begin{theorem}
Let $S$ be a compact oriented surface of genus
$g>1$. Let $G^\C$ be a complex semisimple Lie group and let $G$ be
the split real form of $G^\C$. Then 
the moduli space $\cR(G)$ of representations of the fundamental group
of $S$ in $G$ has a connected  component homeomorphic to a Euclidean
space of dimension $\dim G(2g-2)$. Moreover if $G^\C$ is of adjoint type this
component is unique.
\end{theorem}

When $G=\PSL(2,\R)$, this component can be identified with
Teichm\"uller space, and  this is why Hitchin calls  these
generalized Teichm\"uller components. They are now referred as  
{\bf Hitchin components}.

\subsection{The Hitchin--Kostant--Rallis section}

Hitchin's construction of a section of the Hitchin map for the moduli space 
Higgs bundles for a complex group can be generalized to the moduli space of 
$G$-Higgs bundles for any real reductive Lie group \cite{GPR18}. This construction  
relies  on the Kostant--Rallis construction of a maximal split subalgebra explained 
in Section \ref{kostant-rallis}.

An element $x\in{\liem^\C}$ is said to be \textbf{regular} if 
$\dim \liec_{\liem^\C}(x)=\dim{\liea^\C}$, where 
$\liec_{\liem^\C}(x)=\{y\in\liem^\C\ :\ [y,x]=0\}$.
Denote the subset of regular elements of ${\liem^\C}$ by $\liem^\C_{\reg}.$
Regular elements are those whose ${H^\C}$-orbits are maximal dimensional, so this notion generalises the classical
notion of regularity of an element of a complex reductive Lie algebra.
Note that the intersection $\liem^\C\cap\lieg^\C_{\reg}$ is either empty or the whole of $\liem^\C_{\reg}$. Here 
$\lieg^\C_{\reg}$ denotes the set of elements of $\lieg^\C$ with maximal 
dimensional $G^\C$-orbit. 

A real form $\lieg\subset\lieg^\C$ is called
\textbf{quasi-split} if $\liem^\C\cap\lieg^\C_{reg}=\liem^\C_{reg}$.
 These include
split real forms, and the Lie algebras $\liesu(p,p)$, $\liesu(p,p+1)$, 
$\lieso(p,p+2)$, and $\lie{e}_{6(2)}$. 
Quasi-split real forms admit several equivalent characterizations: $\lieg$ is 
quasi-split if and only if $\liec_\lieg(\liea)$ is abelian --- which holds if and only if $\lieg^\C$ contains a $\theta$-invariant Borel subalgebra --- and 
if and only if $\liem^\C\cap\lieg^\C_{\reg}=\liem_{\reg}$.

Consider the group $Q$, satisfying $(\Ad(G)^\C)^\theta=Q\Ad(H^\C)$.
It is a finite group whose cardinality we denote by $N$. The following theorem is
proved in \cite{GPR18}.
\begin{theorem}\label{}
Let $(G,H,\theta, B)$ be a reductive real Lie group, and let  
$(\widehat{G},\widehat{H},\widehat{\theta},\widehat{B})$
be its maximal split subgroup. Then, the choice of a square root of $K$ 
determines $N$ inequivalent sections of the map
$$
h: \cM(G)\to B(G).
$$
Each such section $s_G$ satisfies
\begin{enumerate}
\item[1.] If $G$ is quasi-split, $s_G(B(G))$ is contained in the stable 
locus of $\cM(G)$, and in the smooth locus
if $Z(G)=Z_G(\lieg)$.
\item[2.] If $G$ is not quasi-split, the image of the section is contained 
in the strictly polystable locus.
\item[3.] For arbitrary groups, the Higgs field is everywhere regular.
\item[4.] The section factors through $\cM(\widehat{G})$. 
\item[5.] If $G_{\Split}<G^\C$ is the split real form satisfying 
$\liez(\lieg^\C)\cap i\lieu\subset\liem$, 
$s_G$ is the factorization of the Hitchin section through $\cM(G_{\Split})$.
\end{enumerate}
\end{theorem}

\section{Hermitian groups and maximal Toledo components}\label{hermitian}

In this section we will assume that $G$ is a connected, non-compact
real simple Lie group of Hermitian type with finite centre.
We fix a maximal compact
subgroup $H\subset G$, with Cartan decomposition $\lieg=\lieh +\liem$.
The centre $Z(H)$ of $H$ is isomorphic
to $\U(1)$. For details on this section  see \cite{BGR17}.

\subsection{Hermitian symmetric spaces and Cayley transform}

The homogenous space $G/H$ is an irreducible Hermitian symmetric space of 
non-compact type. 
The same symmetric space $G/H$ is given by all the finite coverings of
the group $\Ad(G)=G/Z(G)$, where $Z(G)$ is the centre of $G$.  
The vector subspace $\liem$ is isomorphic to the tangent space $T_{o}
G/H$ at the point $o=eH$.
Let $H^\C$, $\liehc$ and $\liemc$ be the complexifications of $H$,
$\lieh$ and $\liem$ respectively.  The almost complex structure $J_0$
on $\liem=T_o(G/H)$, where $o\in G/H$ corresponds to the coset $H$, 
is induced by the adjoint action of an element $J\in
\liez(\lieh)$, so $J_0=\ad(J)|_{\liem}$. Since $J_0^2=-\Id$, we decompose
$\liem^\C$ into $\pm i$-eigenspaces for $J_0$: $\liem^\C=\liem^+
+ \liem^-$. Both $\liem^+$ and $\liem^-$ are abelian, $[\lieh^\C,\liem^\pm]\subset
\liem^\pm$, and there are $\Ad(H)$-equivariant isomorphisms $\liem\cong \liem^\pm$ 
given by $X\mapsto \frac{1}{2}(X\mp iJ_0X)$.

Consider a maximal abelian subalgebra $\liet$ of $\lieh$. Its
complexification $\liet^\C$ gives a Cartan subalgebra of $\lieg^\C$, for
which we consider the root system $\Delta=\Delta(\lieg^\C,\liet^\C)$ and the
decomposition $\lieg^\C=\liet^\C+\sum_{\al\in\Delta} \lieg^\C_\al$. Since
$\ad(\liet^\C)$ preserves $\lieh^\C$ and $\liem^\C$, $\lieg_\al^\C$ must lie
either in $\lieh^\C$ or in $\liem^\C$. If $\lieg_\al^\C\subset \lieh^\C$ (resp.
$\lieg_\al^\C\subset \liem^\C$) we say that the root $\al$ is \textbf{compact}
(resp. \textbf{non-compact}). We choose an ordering of the roots in
such a way that $\liem^+$ (resp. $\liem^-$) is spanned by the root vectors
corresponding to the non-compact positive (resp. negative) roots. 

We denote by $\langle\cdot,\cdot\rangle$ the Killing form on $\lieg^\C$.
For each root $\al\in\Delta$, let $H_\al\in i\liet$ be the dual of $\al$, 
i.e., $$\al(Y)=\langle Y,H_\al\rangle\qquad \textrm{for } Y\in
i\liet.$$ Define, as usual, $h_\al= \frac{2H_\al}{\langle H_\al,H_\al\rangle} \in
i\liet$, and $e_\al\in\lieg^\C_\al$ such that $[e_\al,e_{-\al}]=h_\al$.
Two roots $\al,\be\in\De$ are said to be \textbf{strongly orthogonal}
if neither $\al+\be$ nor $\al-\be$ is a root (equivalently
$[\lieg^{\al},\lieg^{\pm\be}]=\{0\}$). A \textbf{system of strongly orthogonal
  roots} is a maximal set of strongly orthogonal positive non-compact roots.  It
has a number of elements equal to the rank $r=\textrm{rk}(G/H)$ of the
symmetric space $G/H$, i.e., the maximal dimension of a flat, totally
geodesic submanifold of $G/H$. Moreover, for
two strongly orthogonal roots $\ga\neq \ga'$ we have 
\begin{equation}\label{eq:relations-e-h-Gamma}
  [e_{\pm\ga},e_{\pm\ga'}] = 0, \quad [e_{\pm \ga}, h_{\ga'}] = 0.
\end{equation}

For a strongly orthogonal system of roots $\Ga$, consider
\begin{equation*}
  x_\Gamma=\sum_{\ga\in\Gamma} x_\ga, \quad
  y_\Gamma=\sum_{\ga\in\Gamma} y_\ga, \quad
  e_\Gamma=\sum_{\ga\in\Gamma} e_\ga, \quad
  c=\exp\left(\frac{\pi}{4} i y_\Gamma\right) \in U\subset G^\C, 
\end{equation*}
where $G^\C$ is the simply connected Lie group with Lie algebra $\lieg^\C$ and
$U$ is its compact real form  (with Lie algebra $\lieh\oplus i\liem$).  We define the \textbf{Cayley transform} as
the action of the element $c$ on the Lie algebra $\lieg^\C$ by
$\Ad(c):\lieg^\C\to \lieg^\C$. The Cayley transform $\Ad(c)$ satisfies $\Ad(c^8)=\Id$, $\Ad(c)\circ \theta=\theta\circ
\Ad(c^{-1})$ for the Cartan involution $\theta$, and consequently $\Ad(c^4)$ preserves $\lieh$ and $\liem$,
even though $\Ad(c)$ does not preserve $\lieg$.

It is well-known that a Hermitian symmetric space of non-compact type
$G/H$ can be realized as a bounded symmetric domain.
For the classical groups this
is due to Cartan, while the  general case
is given by the Harish-Chandra embedding $ G/H \to \liemp$
which defines
a biholomorphism between $G/H$ and
the bounded symmetric domain $\cD$ given  by the image of $G/H$ in
the complex  vector space $\liemp$.
Now, for any bounded domain $\cD$ there is  the {\bf Shilov boundary}
of $\cD$ which is defined as the smallest closed subset $\Sh$
of the topological boundary $\partial \cD$ for which  every function
 $f$ continuous on $\overline{\cD}$ and holomorphic on $\cD$ satisfies
that
$$
|f(z)|\leq \max_{w \in \Sh} |f(w)|\;\;\mbox{for every}\;\; z\in\cD.
$$
The Shilov boundary $\Sh$ is the unique closed $G$-orbit in
$\partial \cD$.

The simplest situation to consider is that of the hyperbolic plane.
The Poincar\'e disc is its realization as a bounded symmetric domain.
However, we know that the hyperbolic plane  can also be realized as
the Poincar\'e upper-half plane. There are other Hermitian symmetric spaces that,
like the hyperbolic plane, admit a realization similar to   the
upper-half plane. These are the tube type symmetric spaces.

Let $V$ be a real vector space and let $\Omega \subset V$ be an open cone
in $V$. A tube over the cone $\Omega$  is a domain of the form
$$
T_\Omega=\{u+iv\in V^\C, u\in V, v\in \Omega\}.
$$
A domain $\cD$ is said to be of {\bf tube type}  if it is biholomorphic
to a tube $T_\Omega$. In the case of a symmetric domain the cone $\Omega$
is also symmetric.
An important characterization of the tube type symmetric domains is given
by the following.

\begin{proposition}Let $\cD$ be a bounded symmetric domain corresponding to 
the Hermitian symmetric space $G/H$. The following are equivalent:

(i) $\cD$ is of tube type.

(ii) $\dim_\R \Sh=\dim_\C\cD$.

(iii) $\Sh$ is a symmetric space of compact type.

(iv) $\Ad(c^4)=\Id$.
\end{proposition}

There is a generalization of the Cayley map that sends the  unit disc
biholomorphically to the upper-half plane. Let  $\cD$ be  the bounded domain
associated to a Hermitian  symmetric space $G/H$.
Acting by  a particular element in $G$, known as the {\bf Cayley element},
one obtains a map
$$
c: \cD\lra \liem^+
$$
which is called  the  {\bf Cayley transform}.
A  relevant fact for us is the following.

\begin{proposition}\label{tube-domain}
Let $\cD$ be the  symmetric domain corresponding to the Hermitian
symmetric space $G/H$. 
Let $\cD$ be of tube type. Then the image by
the Cayley transform $c(\cD)$  is biholomorphic to a tube domain
$T_\Omega$ where the symmetric cone $\Omega$ is the non-compact
dual of the  Shilov boundary of $\cD$. In fact the Shilov boundary
is a symmetric space isomorphic to  $H/H'$ for a certain subgroup
$H'\subset H$, and $\Omega=H^*/H'$ is its non-compact dual symmetric
space.
\end{proposition}

\begin{proposition}
  \begin{enumerate}
  \item The symmetric spaces defined by $\Sp(2n,\R)$, $\SO_0(2,n)$ are
    of tube type.
  \item The symmetric space defined by $\SU(p,q)$ is of tube type if
    and only if $p=q$.
  \item The symmetric space defined by $\SO^*(2n)$ is of tube type if
    and only if $n$ is even.
  \end{enumerate}
\end{proposition}

For a  tube type classical irreducible symmetric space $G/H$,
Table~\ref{tab:tube} indicates  the Shilov boundary
$\Sh=H/H'$, its non-compact dual $\Omega= H^*/H'$,
the isotropy representation space $\liem'$
and its complexification ${\liem'}^\C$, corresponding to the Cartan
decomposition of the Lie algebra $\Lie(H^*)=\lieh'+\liem'$ of $H^*$.
The vector space $\liem'$  has the structure of a Euclidean Jordan
algebra,  where the cone $\Omega$ is realized.

The following lemma will be important for our Cayley correspondence.

\begin{lemma}\label{lemma:cayley-iso}
  There are $\Ad(H'^\C)$-equivariant isomorphisms $\ad(e_\Gamma):\liem_T^-
  \to \liem'^\C$ and $\ad(e_\Gamma):\liem'^\C \to \liem_T^+$. 
\end{lemma}

The study of certain problems in non-tube type domains
can be reduced to the tube type thanks to the following.

\begin{proposition}
Let $G/H$ be  a Hermitian symmetric space of non-compact type.
There exists a subgroup  $\widetilde G\subset G$ such that
$\widetilde G/\widetilde H\subset
G/H$ is a maximal isometrically embedded  symmetric space of tube type,
where  $\widetilde H\subset \widetilde G$ is a maximal compact subgroup.
\end{proposition}

Table~\ref{tab:non-tube} gives the maximal symmetric space of tube
type isometrically embedded in the two series of irreducible classical
symmetric spaces of non-tube type. We describe also the Shilov
boundaries of $G/H$ and $\widetilde G/\widetilde H$ which are of the form
$\Sh=H/H'$, and $\widetilde{\Sh}=\widetilde H/\widetilde H'$, respectively. Notice
that in the non-tube case the Shilov boundary $\Sh$ is a homogeneous
space $H/H'$ but it is not symmetric.

\subsection{The Toledo character}

We introduce the Toledo character associated to a simple Lie algebra $\lieg$ of Hermitian type as the character on the Lie algebra $\lieh^\C$ given as follows.

\begin{definition}
  The \textbf{Toledo character} $\chi_T:\lieh^\C \to \C$ is defined, for $Y\in \lieh^\C$ in terms of the Killing form, by
$$ \chi_T(Y) = \frac1N \langle-i J, Y\rangle, $$
where $N$ is the dual Coxeter number. 
\end{definition}
Since $J$ is in the center, $\chi_T$ vanishes on $[\lieh^\C,\lieh^\C]$, hence determines a character.

We study now when the Toledo character lifts to a character of the
group $H^\C$. Note that this depends on the choice of the pair $(G,H)$
defining the same symmetric space. Let $Z_0^\C$ denote the identity
component of $Z(H^\C)$.
\begin{proposition}\label{prop:exponentiation-of-Toledo}
Define $o_J$ to be the order of $e^{2\pi J}$ and $\ell=|Z_0^\C\cap
[H^\C,H^\C]|$. For $q\in \Q$, the character $q \chi_T$ lifts to $H^\C$ if
and only if $q$ is an integral multiple of
\begin{equation}\label{eq:Toledo-character-group}
q_T=\frac{\ell N}{o_J \dim \liem}.
\end{equation}
\end{proposition}

The value of $q_T$ in the standard examples is given in Table 
\ref{tab:exponent-tubetype}. 
Note in particular that $q_T=\frac12$ for
all classical groups except $\SO^*$, for which $q_T=1$. So for all
classical groups the Toledo character lifts to $H^\C$.
For the adjoint group, $o_J=1$ so $q_T=\frac{\ell N}{\dim\fm}$. In the tube case $N=\frac{\dim\fm}r$ so this gives $q_T=\frac \ell r$. The values of 
$q_T$ in the non-tube case are given in Table \ref{tab:exponent-adjoint}.

Finally, the following lemma will prove later that the Toledo invariants 
defined from the two points of view of Higgs bundles and representations 
coincide.
\begin{lemma}\label{lemma:Toledo-Kahler-form}
 The Toledo character $\chi_T$ defines a symmetric K\"ahler form on $G/H$ by 
$$\omega(Y,Z)=i\chi_T([Y,Z]), \text{ for }Y,Z\in \liem,$$
with minimal  holomorphic sectional curvature $-1$. 
\end{lemma}

We define now for $G$ of tube type a determinant polynomial, $\det$, on 
$\liem^+$, whose degree
equals the rank of the symmetric space. This determinant is a familiar object
in Jordan algebra theory \cite{FK94}, but it can be introduced in an
elementary way as follows \cite[Lemma 2.3]{KV79}: it is the unique
$H'_0$-invariant 
polynomial on $\liem^+$ which restricts on $\liea^+=\oplus_1^r\C e_{\gamma_i}$ to
$$ \det \sum_1^r \lambda_ie_{\gamma_i} = \prod_1^r \lambda_i. $$
Here $H'_0$ is the identity component of $H'$.
The existence comes from the Chevalley theorem on invariant polynomials, since the Weyl group acts exactly by all permutations on the $(e_{\gamma_i})$ (see again \cite{KV79}).

The main useful property for us is the following equivariance:
\begin{lemma}\label{lemma:relation-det-char-to-q}
 Let $G$ be of tube type. For $h\in H^\C$ and $x\in \liem^+$ we have
\begin{equation}\label{eq:8}
  \det(\Ad(h) x)=\widetilde{\chi}_T(h) \det(x),
\end{equation}
where $\widetilde{\chi}_T$ is the lifting of $\chi_T$ to $H^\C$.
\end{lemma}
Note that we implicitly assumed here that the lifting $\widetilde{\chi}_T$
exists, otherwise the same identity remains true after taking power
$q_T$. 

We define a notion of rank for an element in  $\liem^+$
for $G$ of Hermitian type (tube or non-tube). 
Choose $\liea^+=\oplus_1^r\C e_{\gamma_i}$. Any element of $\liem^+$ 
is conjugate under $H^\C$ to an element of $\liea^+$.
\begin{definition}\label{definition-rank}
Let $x\in \liem^+$, and $y=\sum \lambda_ie_{\gamma_i}\in \liea^+$ be conjugate to $x$ under $H^\C$. Then we say that $x$ has rank $r'$ if $y$ has exactly $r'$ non-zero coefficients.
\end{definition}
This is well defined because the Weyl group acts only by permutations
on $\liea^+$. Also, in the tube case, one can give a more
intrinsic interpretation using the determinant: polarize the
determinant to get an $r$-linear map $C$ on $\liem^+$ such that
$C(x,\ldots,x)=\det(x)$; then the rank of $x$ is the maximal 
integer $r'$ such that the $(r-r')$-form $C(x,\ldots,x,\cdot,\ldots,\cdot)$ 
is not identically zero,
which is clearly an invariant notion.

\begin{remark}
 In the case of $\SU(p,q)$ the rank on $\liem^+$ specializes to the notion of rank for a rectangular matrix $q\ts p$. For $\Sp(2n,\R)$, the rank on $\liem^+$ is the rank for an element of $S^2(\C^n)$ seen as an endomorphism. 
\end{remark}

The following proposition plays an important role in what follows.
\begin{proposition}\label{prop:HC-trans}
Let $1\leq r'\leq r$.  The group $H^\C$ acts transitively on the set of  elements of rank $r'$ in $\liem^+$. In particular, the set of regular (that is maximal rank) elements in $\liem^+$ is $H^\C/H'^\C$.
\end{proposition}

\subsection{Toledo invariant and Milnor--Wood inequality}

Let $X$ be a compact Riemann surface and let $(E,\varphi)$ be a $G$-Higgs 
bundle over $X$.  The decomposition
$\liem^\C=\liem^++\liem^-$ gives a vector bundle decomposition $E(\liem^\C)=
E(\liem^+) \oplus E(\liem^-)$ and hence the Higgs field has two components:
$$
\varphi=(\varphi^+, \varphi^-)\in H^0(X,E(\liem^+)\otimes K)\oplus
H^0(X,E(\liem^-)\otimes K)= H^0(X,E(\liem^\C)\otimes K).
$$

When the group $G$ is a classical group, or more generally when
$H$ is a classical group, 
it is  useful to take the standard representation of $H^\C$ to describe a 
$G$-Higgs bundle in terms of associated vector bundles. 
This is the approach taken in \cite{BGG03,BGG06,BGG15,GGM13}. 

Let $(E,\varphi)$ be a $G$-Higgs bundle. Consider the Toledo character $\chi_T$. Up to an integer
multiple, $\chi_T$ lifts to a character $\tilde{\chi}_T$ of $H^\C$. Let
$E(\tilde{\chi}_T)$ be the line bundle associated to $E$ via the
character $\tilde \chi_T$.

\begin{definition}\label{def:Toledo-invariant}
  We define the \textbf{Toledo invariant} $\tau$ of $(E,\varphi)$ 
by\newnot{Toledoinvariant}
$$
\tau=\tau(E):=\deg(E(\tilde{\chi}_T)).
$$
\end{definition}

If $\tilde \chi_T$ is not defined, but only $\tilde \chi_T^q$, one must replace the definition by $\frac1q\deg E(\tilde{\chi}^q_T).$

We denote by $\cM^\alpha_\tau(G)$ the subspace of $\cM^\alpha(G)$ 
corresponding to  $G$-Higgs bundles whose Toledo invariant equals $\tau$.
For $\alpha=0$ we simplify our notation setting  $\cM_\tau(G):= \cM^0_\tau(G)$

The following proposition relates our Toledo invariant to the usual Toledo invariant of a representation, first defined in \cite{toledo}.
\begin{proposition}\label{prop:toledo-rep}
  Let $\rho:\pi_1(X) \to G$ be reductive and let $(E,\varphi)$ be the corresponding
  polystable $G$-Higgs bundle given by Theorem \ref{na-Hodge}. Let $f:\tilde X\to G/H$ be the corresponding harmonic metric. Then
$$ \tau(E) = \frac1{2\pi} \int_X f^*\omega , $$
  where $\omega$ is the K\"ahler form of the symmetric metric on $G/H$ with
  minimal holomorphic sectional curvature $-1$, computed
 in Lemma \ref{lemma:Toledo-Kahler-form}.

  In particular, $\tau(E)$ is the Toledo invariant of $\rho$.
\end{proposition}

 The Toledo invariant  is related to the topological class of the
 bundle $E$ defined as an element of $\pi_1(H)$. To explain this, assume that
 $H^\C$ is connected. The topological classification of $H^\C$-bundles $E$ on
 $X$  is given by a characteristic class $c(E)\in \pi_1(H^\C)$ as follows. 
 From the exact sequence 
 $$
 1\to \pi_1(H^\C) \to \widetilde{H^\C} \to H^\C \to 1
 $$  
 we obtain a long
 exact  sequence in cohomology and, in particular, the connection map 
 \begin{equation}
 H^1(X,\underline{H}^\C)\xra{c} H^2(X,\pi_1(H^\C)),\label{eq:1}
 \end{equation}
 where $\underline{H}^\C$ is the sheaf of local holomorphic functions in $X$  
with 
 values in $H^\C$. The cohomology set $H^1(X,\underline{H}^\C)$ 
 (not necessarily  a group since $H^\C$ is  in general
 not abelian) parametrizes isomorphism classes of principal $H^\C$-bundles 
 over $X$. On the other hand, since $\dim_\R X=2$, by the universal 
coefficient theorem  and the fact that the fundamental group of a 
Lie group is abelian,    $H^2(X,\pi_1(H^\C))$ is isomorphic to
 $\pi_1(H^\C)$. Moreover, $\pi_1(H^\C)\cong \pi_1(H)\cong\pi_1(G)$ 
since $H$ is a deformation retract for
 both $H^\C$ and $G$. This map thus associates a topological invariant 
in $\pi_1(H)$  to any $G$-Higgs bundle on $X$.

 By the relation between the fundamental group and the centre of a Lie
 group, the topological class in $\pi_1(H)$ is of special interest when
 $H$ has a non-discrete centre, i.e., when $G$ is of Hermitian type. In
 this case, $\pi_1(H)$ is isomorphic to $\Z$ plus possibly a torsion
 group (among the classical groups, $\SO_0(2,n)$ is the only one with
 torsion). Very often (see for example \cite{BGG06}), the Toledo
 invariant of a $G$-Higgs bundle $(E,\varphi)$  is defined as the projection
 of $c(E)$ defined by (\ref{eq:1}) on the torsion-free part, $\Z$. The
 general relation is the following.

 \begin{proposition}\label{prop:Milnor-Wood-mod-invariant}
  Let $(E,\varphi)$ be a $G$-Higgs bundle, and $d\in \Z$ the projection on the
  torsion-free part of the class $c(E)$ defined by (\ref{eq:1}). Then
  $d$ is related to the Toledo invariant by
  $$ 
 \tau = \frac{d}{q_T}.$$
  \end{proposition}

Definition \ref{definition-rank} gives the ranks of $\varphi^+$  and $\varphi^-$ 
at a point $x\in X$. The space of elements of $\liem^\pm$ with rank at most $\rho$ 
is an algebraic subvariety of $\liem^\pm$, so the ranks of $\varphi^+$ and $\varphi^-$ 
are the same at all points of $X$ except a finite number of points where it it smaller. 
We therefore have a well defined notion of rank of $\varphi^+$ and $\varphi^-$:
\begin{definition}\label{def:rank-Higgs-field}
  The generic value on $X$ of the rank of $\varphi^+$ is called the rank of
  $\varphi^+$ and denoted $\rk \varphi^+$. Analogously we define the rank $\rk \varphi^-$ of
  $\varphi^-$.
\end{definition}

The main result of this section is the following.
\begin{theorem}\label{0-theo:ineq-rk}   Let $(E,\varphi^+,\varphi^-)$ be a
semistable $G$-Higgs  bundle. Then, the Toledo invariant of $E$ satisfies:
$$
-\rk(\varphi^+)(2g-2) \leq
\tau \leq \rk(\varphi^-)(2g-2).
$$
In particular, we obtain the familiar Milnor--Wood inequality
$$
|\tau| \leq \rk(G/H)(2g-2),
$$
and  equality holds if and only if $\varphi^+$ 
(resp. $\varphi^-$) is regular at each point in the case $\tau<0$ 
(resp. $\tau>0$).
\end{theorem}

Theorem \ref{0-theo:ineq-rk} 
was proved on a case by case 
basis for the classical groups \cite{hitchin87,gothen,BGG03,BGG06,BGG15,GGM13}. 
In these 
references, the bound  given is for the integer 
$d\in \pi_1(H)\cong\pi_1(H^\C)\cong\Z$ associated 
naturally to the $H^\C$-bundle $E$. This differs
from the Toledo invariant by a rational multiple.
From Table \ref{tab:exponent-tubetype} 
and Proposition \ref{prop:Milnor-Wood-mod-invariant} combined with 
Theorem \ref{0-theo:ineq-rk} we obtain  the Milnor--Wood  inequalities given
 in \cite{BGG06} for the classical Hermitian groups.
Our intrinsic general approach covers of course the exceptional
groups and quotients and covers of classical groups that have not
been studied previously. 

A polystable $G$-Higgs bundle
$(E,\varphi)$ is, by  Theorem \ref{na-Hodge},  in correspondence with a reductive representation
$\rho:\pi_1(X)\to G$, and from Proposition  \ref{prop:toledo-rep}
the Toledo invariant of $(E,\varphi)$ coincides
with the Toledo invariant of a representation of the fundamental 
group in $G$. 
In the context of representations the inequality $|\tau|\leq \rk(G/H)(2g-2)$,
goes back to Milnor \cite{milnor}, who studies the case $G=\PSL(2,\R)$,  and was
proved  in various cases in  \cite{wood,dupont,DT87,CO03}, and in general 
in \cite{BIW10}.
We should point out that the Higgs bundle 
approach gives the Milnor--Wood inequality for an arbitrary representation, 
as the other approaches do,  since such a representation   
can always be deformed to a reductive one. 

\subsection{Hermitian groups of tube type and Cayley correspondence}
\label{sec:cayley}
We define
 a polystable Higgs bundle $(E,\varphi)$ to be \textbf{maximal} if its
 Toledo invariant $\tau$ attains one of the bounds of the inequality i.e., $\tau=\pm
 r(2g-2)$, where $r=\rk(G/H)$. We denote $\tau_{\max}=\rk(G/H)(2g-2)$.

Let $H^*$ be the non-compact dual of $H$ as defined in Proposition 
\ref{tube-domain}.  In this section we establish a bijective
correspondence between maximal $G$-Higgs bundles over $X$ and
$K^2$-twisted $H^*$-Higgs bundles over $X$, as defined in Remark
\ref{twisted}, where $K^2$ is the square of the canonical line bundle.

Suppose that $(E,\varphi)$ is a polystable maximal $G$-Higgs bundle, and choose for example $\tau=-r(2g-2)$. By Theorem \ref{0-theo:ineq-rk}, the field $\varphi^+$ has rank $r$ at each point. Let $Z_0^\C\simeq \C^*$ be the connected component of the identity of the center of $H^\C$. There is an exact sequence
\begin{equation}
 1 \longrightarrow Z_0^\C \longrightarrow H^\C \longrightarrow H^\C/Z_0^\C \longrightarrow 1,\label{eq:13}
\end{equation}
so there is an action of $Z_0^\C$-bundles on $H^\C$-bundles that we will
denote by $\otimes$: in this way, if $\kappa$ is a line bundle over $X$, we
can 
define $E\otimes \kappa$ (here we are identifying the line bundle $\kappa$
with its corresponding $\C^\ast$-bundle).

  If $\kappa$ is an $o_J$-root of $K$, where $o_J$ is the order of $e^{2\pi
    J}$, then $\varphi^+$ defines a reduction of the $H^\C$-bundle $E\otimes
  \kappa$ 
to the group $H'^\C$.

Of course, such $\kappa$ exists only if $o_J$ divides $2g-2$. We will now suppose that $\kappa$ exists and is fixed. Denote by $E'$ the reduction of $E\otimes \kappa$ to $H'^\C$. As we have seen, $\varphi^+\in H^0(X,E'(\fm^+))$, and similarly $\varphi^-\in H^0(X,E'(\fm^-)\otimes K^2)$. 
From Lemma \ref{lemma:cayley-iso}, we have an isomorphism
\begin{equation}
  \label{eq:11}
  \ad \varphi^+ : E'(\fm^-) \longrightarrow E'(\fm'^\C),
\end{equation}
so that we can define a Higgs field 
$$\varphi'=[\varphi^+,\varphi^-]\in H^0(X,E'(\fm'^\C)\otimes K^2).$$ The data $(E',\varphi')$ is a $K^2$-twisted $H^*$-Higgs bundle.

Conversely, from a $K^2$-twisted $H^*$-Higgs bundle $(E',\varphi')$ we can reconstruct $(E,\varphi)$ in the following way. The bundle is $E=E'\otimes \kappa^{-1}$. Observe that for the $H'^\C$-bundle $E'$ we have canonical section $e_\Gamma\in H^0(X,E'(\fm^+))$ corresponding to the element $e_\Gamma\in \fm^+$ fixed by $H'$, which becomes by Lemma \ref{lemma:cayley-iso} a section $\varphi^+\in H^0(X,E(\fm^+)\otimes K)$. Finally, $\varphi^-$ is reconstructed from (\ref{eq:11}) as $(\ad \varphi^+)^{-1}(\varphi')$. Therefore, $\kappa$ being fixed, we obtain a complete correspondence between maximal $G$-Higgs bundles and $K^2$-twisted $H^*$-Higgs bundles.
We refer to  $(E',\varphi')$ as the \textbf{Cayley partner} of $(E,\varphi)$.
One hs the following.

\begin{theorem}[{\bf Cayley correspondence}]  
\label{th:cayley-correspondence}
Let $G$ be a connected non-compact  real simple Hermitian Lie group of tube type
with finite centre. Let $H$ be a maximal compact subgroup of $G$ and  $H^*$ be 
the non-compact dual of $H$ in $H^\C$. 
Let $J$ be the element in $\liez$ (the centre of $\lieh$) defining
the almost complex structure on $\liem$. 
If the order of $e^{2\pi J}\in H^\C$ divides $(2g-2)$, then there is an 
isomorphism of complex algebraic varieties
\begin{equation}\label{eq:generalized-cayley}
\cM_{\max} (G) \cong  \cM_{K^2}(H^*)
\end{equation}
given by $(E,\varphi)\mapsto(E',\varphi')$ as above.
\end{theorem}

\begin{remark}
  The condition $o_J | (2g-2)$ is always satisfied for a group of
  adjoint type, since in this case $o_J=1$.  Table
  \ref{tab:exponent-tubetype}  shows
  that the $o_J$ divides $(2g-2)$ for the classical and exceptional
  groups. This  may not happen for coverings of these groups, where
  $o_J$ may be bigger.
\end{remark}

\section{$\SO(p,q)$-Higgs bundles and higher Teichm\"uller spaces} 
Details on this section can be found in \cite{A-et-al1,A-et-al2}.
\subsection{$\SO(p,q)$-Higgs bundles and topological invariants} 

An $\SO(p,q)$-Higgs bundle on $X$ is equivalent to a  triple $(V,W,\eta)$ where 
$V$ and $W$ are respectively rank $p$ and rank $q$ vector bundles with orthogonal structures such that $\det(W)\simeq\det(V)$, and $\eta$ is a holomorphic bundle map $\eta:W\rightarrow V\otimes K$.  

The cases for $p\leq2$ are somewhat special. For $p>2,$ rank $p$ orthogonal bundles on $X$ are classified  topologically by their first and second Stiefel-Whitney classes, $sw_1\in H^1(X,\Z_2)$ and $sw_2\in H^2(X,\Z_2)$. These primary topological invariants are constant on connected components of the moduli space $\mathcal{M}(\SO(p,q))$.  Since $\det(W)\simeq\det(V)$, it follows that $sw_1(V)=sw_1(W)$.  The components of the moduli space $\mathcal{M}(\SO(p,q))$ are thus partially labeled by triples $(a,b,c)\in \Z_2^{2g}\times\Z_2\times\Z_2$, where 
\[\xymatrix@C=.3em{a=sw_1(V)\in H^1(X,\Z_2),&b=sw_2(V)\in H^2(X,\Z_2)&\text{and}&c=sw_2(W)\in H^2(X,\Z_2).}\]

Using the notation $\mathcal{M}^{a,b,c}(\SO(p,q))$ to denote the union of components labeled by $(a,b,c)$, we can thus write
\begin{equation}\label{Mpq-abc}
\mathcal{M}(\SO(p,q))=\coprod_{(a,b,c)\in \Z_2^{2g}\times\Z_2\times\Z_2}\mathcal{M}^{a,b,c}(\SO(p,q))\ .
\end{equation}

Stability  for $\SO(p,q)$-Higgs bundles implies that a Higgs bundle $(V,W,\eta)$ with $\eta=0$ is polystable if and only if $V$ and $W$ are both polystable orthogonal bundles. This leads to the immediate identification of one connected component in each space $\mathcal{M}^{a,b,c}(\SO(p,q))$.  
We use the subscript `top' to designate these components, which contain $\SO(p,q)$-Higgs bundles with vanishing Higgs field.

\begin{proposition}\label{zerocomponents}  Assume that $2<p\le q$.  
For every $(a,b,c)\in \Z_2^{2g}\times\Z_2\times\Z_2$ the space $\mathcal{M}^{a,b,c}(\SO(p,q))$ has a non-empty connected component,
denoted by $\mathcal{M}_{\rmin}^{a,b,c}(\SO(p,q))$, in which every point can be continuously deformed to the isomorphism class of an $\SO(p,q)$-Higgs bundle of the form $(V,W,\eta=0)$ where $V$ and $W$ are polystable orthogonal bundles.
\end{proposition}
We define
\begin{equation}
\mathcal{M}_{\rmin}(\SO(p,q))=\coprod_{a,b,c}  \mathcal{M}_{\rmin}^{a,b,c}(\SO(p,q))
\end{equation}

\begin{remark}
In the case $p=2$ it is no longer true that $\mathcal{M}_{\rmin}^{a,b,c}(\SO(p,q))$ is non-empty for all $(a,b,c)$. In particular, if $a=0,$ then $V=L\oplus L^{-1}$ which (a) is polystable if $\deg{L}=0$ and (b) has $sw_2(V)=\deg{L} \mod 2$. Thus $\mathcal{M}_{\rmin}^{0,b,c}\SO(2,q))$ is empty if $b\ne 0$.
\end{remark}

The  main result in \cite{A-et-al1,A-et-al2} is  that the moduli space $\mathcal{M}(\SO(p,q))$ has additional `exotic' components 
disjoint from the components 
of $\mathcal{M}_{\rmin}(\SO(p,q))$.  These exotic components are identified as products of moduli spaces of so-called $L$-twisted Higgs bundles, where in each factor $L$ is a positive power of the canonical bundle $K$.  

Let $L$ be a fixed holomorphic line bundle on $X$.  An $L$-twisted 
$\SO(1,n)$-Higgs bundle on $X$ is equivalent to a triple $(I,W_0,\eta)$, where $W_0$ is a $\Or(n,\C)$-bundle, $I$ is the rank one orthogonal bundle $\det(W_0)$ and $\eta:W_0\to I\otimes L$ is a holomorphic bundle map. 

We get a decomposition similar to \eqref{Mpq-abc}, namely
\begin{equation}\label{M1n-ac}
\mathcal{M}_{K^p}(\SO(1,n))=\coprod_{(a,c)\in \Z_2^{2g}\times\Z_2}\mathcal{M}_{K^p}^{a,c}(\SO(1,n))\ ,
\end{equation}
\noi where $\mathcal{M}_{K^p}^{a,c}(\SO(1,n))$ denotes the component in which the $\SO(1,n)$-Higgs bundles are of the form $(I,W_0,\eta)$, with $a=sw_1(W_0)$ and $c=sw_2(W_0)$.

\subsection{The case $2<p\le q$}  

We can now state the  main result in \cite{A-et-al1,A-et-al2}.

\begin{theorem}\label{mainth} Fix integers $(p,q)$ such that $2< p< q-1$.  For each choice of $a\in \Z_2^{2g}$ and $c\in\Z_2$, the moduli space $\mathcal{M}(\SO(p,q))$ has a connected component disjoint from $\mathcal{M}_{\rmin}(\SO(p,q))$. This component is isomorphic to 
\begin{equation}\label{pqcayley}
\mathcal{M}_{K^p}^{a,c}(\SO(1,q-p+1))\times\mathcal{M}_{K^2}(\SO_0(1,1))\times\cdots\times\mathcal{M}_{K^{2p-2}}(\SO_0(1,1))~,
\end{equation}
and lies in the sector $\mathcal{M}^{\alpha,0,c}(\SO(p,q))$ where $\alpha=a$ if $p$ is odd and $\alpha=0$ if $p$ is even. Moreover, $\cM(\SO(p,q))$ has no other connected components. 
 \end{theorem}

\begin{remark}
The group $\SO_0(1,1)$ is the connected component of the identity in $\SO(1,1)$, and $\mathcal{M}_{K^{2j}}(\SO_0(1,1))$ can be identified with $H^0(K^{2j})$. Thus, we can replace \eqref{pqcayley} with
\begin{equation}\label{pqcayley.2}
\mathcal{M}_{K^p}^{a,c}(\SO(1,q-p+1))\times \bigoplus_{j=1}^{p-1}H^0(K^{2j})~.
\end{equation}
\end{remark}

\begin{remark}
The existence of the exotic components described by \eqref{pqcayley} was 
proven for $p=2$ in \cite{BGG06}. They are the exotic components with maximal Toledo invariant arising from Cayley correspondence (see Section \ref{SO2q}). 
In particular, Theorem \ref{mainth} can be viewed as a generalized Cayley correspondence. Contrary to the cases $p>2$, there are components of $\mathcal{M}(\SO(2,q))$ which are not in the family described by the theorem and also not in $\mathcal{M}_{\rmin}(\SO(2,q))$. These are the components with non-maximal and non-zero Toledo invariant. 
\end{remark}

\begin{corollary} For $2<p<q-1$, the moduli space $\mathcal{M}(\SO(p,q))$ has $3\times 2^{2g+1}$ connected components, $2^{2g+1}$ of which are exotic components disjoint from $\mathcal{M}_\mathrm{top}(\SO(p,q))$. 
\end{corollary}

\subsection{The case $q=p+1$}  \label{brian}

If $q=p+1$, then $\mathcal{M}_{K^p}^{a,c}(\SO(1,q-p+1))=\mathcal{M}_{K^p}^{a,c}(\SO(1,2))$, which is not always 
connected.  Indeed, if $a=0$, then the Higgs bundles represented in $\mathcal{M}_{K^p}^{0,c}(\SO(1,2))$ can be taken to be of the form $(\mathcal{O}, L\oplus L^{-1}, \eta)$, where $L$ is a non-negative degree $d$ line bundle. Stability considerations impose a bound on $d$ so that 
\begin{equation}\label{so120}
\mathcal{M}_{K^p}^{0,c}(\SO(1,2))=\coprod_{\substack{0\le d\le p(2g-2)\\ d=c\,(\mathrm{mod}\ 2)}}\mathcal{M}_{K^p}^{d}(\SO(1,2)).
\end{equation}
Moreover, Collier \cite{Col17} has showon that  for each integer $d\in (0, 2g-2]$, the moduli space $\mathcal{M}_{K^p}^{d}(\SO(1,2))$ is diffeomorphic to a vector bundle of rank $d+g-1$ over the $(2g-2-d)^{th}$-symmetric product $\Sym^{2g-2+d}(X)$. In particular, the components $\mathcal{M}_{K^p}^{d}(\SO(1,2))$ are smooth and connected.  

The moduli spaces $\mathcal{M}(\SO(p,p+1))$ have been analyzed by Collier  
\cite{Col17}. It was shown there that the topological invariants for $\SO(p,p+1)$-Higgs bundles, i.e.\  the triples $(a,b,c)$, do not distinguish all connected components.  Two families of exotic components were identified. The components in the first family are labeled by an integer, $d$, in the range $0\le d\le p(2g-2)$, while those in the second family are labeled by a pair $(a,c)\in(\Z^{2g}-\{0\})\times\Z_2$. Though not described in this way in \cite{Col17}, these families can be identified as follows:

\begin{itemize}
\item In the the family labeled by $d$, each member is isomorphic to
\begin{equation}
\mathcal{M}_{K^p}^{d}(\SO(1,2))\times\mathcal{M}_{K^2}(\SO_0(1,1))\times\cdots\times\mathcal{M}_{K^{2p-2}}(\SO_0(1,1)),
\end{equation}
\noi where $\mathcal{M}_{K^p}^{d}(\SO(1,2))$ is one of the components of $\mathcal{M}_{K^p}^{0,c}(\SO(1,2))$ as in \eqref{so120}.

\item In the family labeled by $(a,c)$, each member is isomorphic to
\begin{equation}
\mathcal{M}_{K^p}^{a,c}(\SO(1,2))\times\mathcal{M}_{K^2}(\SO_0(1,1))\times\cdots\times \mathcal{M}_{K^{2p-2}}(\SO_0(1,1)).
\end{equation}
\end{itemize}

\noi The components are thus precisely those identified by Theorem \ref{mainth} in the case  $q=p+1$.  The component count in this case is, however, different from the case $q>p+1$.

\begin{corollary} For $p>2$, the moduli space $\mathcal{M}(\SO(p,p+1))$ has $3\times 2^{2g+1}+2p(g-1)-1$ connected components.  Among those, there are $2^{2g+1}+2p(g-1)-1$ `exotic' components which are disjoint from $\cM_{\rmin}(\SO(p,p+1))$.  
\end{corollary}

\subsection{The case $q=p$}

 In this case $\mathcal{M}_{K^p}(\SO(1,q-p+1))=\mathcal{M}_{K^p}(\SO(1,1))$.  A $K^p$-twisted $\SO(1,1)$-Higgs bundle consists of a triple $(I,I,\eta)$ where $I$ is a square root of the trivial bundle $\mathcal{O}$ and $\eta\in H^0(K^p)$.  Such Higgs bundles are labeled by a single Stiefel-Whitney class, namely $a=sw_1(I)$, so that

\begin{equation}
\mathcal{M}_{K^p}(\SO(1,1))=\coprod_{a\in H^1(X,\Z_2)}\mathcal{M}^a_{K^p}(\SO(1,1)).
\end{equation}

\noi With $q=p,$ Theorem \ref{mainth} thus gives $2^{2g}$ exotic components of $\mathcal{M}(\SO(p,p))$ isomorphic to the moduli spaces
\begin{equation}
\mathcal{M}_{K^p}^{a}(\SO(1,1))\times\mathcal{M}_{K^2}(\SO_0(1,1))\times\cdots\times\mathcal{M}_{K^{2p-2}}(\SO_0(1,1)).
\end{equation}

\noi For each $a$, we can identify $\mathcal{M}_{K^p}^{a}(\SO(1,1))$ with $H^0(K^p)$. Thus, each exotic component is isomorphic to $H^0(K^p)\oplus\bigoplus\limits_{j=1}^{p-1}(H^0(K^{2j})$.  This recovers 
 the Hitchin component in $\mathcal{M}(\SO_0(p,p))$ when $a=0$.

\subsection{The case $p=2<q$.}\label{SO2q} 

An $\SO(2,q)$-Higgs bundle is defined by a triple $(V,W,\eta)$ in which $V$ is an $\mathrm{O}(2,\C)$-bundle.  If $sw_1(V)=0$, i.e.\ if the structure group of $V$ reduces to $\SO(2,\C)$, then $V$ can be assumed to be a direct sum of line bundles of the form $V=L\oplus L^{-1}$, with orthogonal structure given $q_V=\smtrx{0&1\\1&0}$ in this splitting. Note that the second Stiefel-Whitney class of the orthogonal bundle $L\oplus L^{-1}$ is given by $sw_2=d \,(\mathrm{mod}\ 2)$ where $d=\deg(L)\geq0$.  

For the groups $\SO(2,q)$, the connected components of the identity are isometry groups of Hermitian symmetric spaces of non-compact type. As explined in Section \ref{hermitian}, the Higgs bundles  have an associated Toledo invariant which, up to a normalization constant, is integer-valued but subject to a Milnor-Wood bound. For an $\SO_0(2,q)$-Higgs bundle $(L\oplus L^{-1}, W, \eta)$, the Toledo invariant is basically the degree $d$ of $L$ and the Milnor-Wood bound is $0\leq d\le 2g-2$.
We thus get
\begin{equation}
\mathcal{M}^{0,b,c}(\SO(2,q))=\coprod_{\substack{0\leq d\leq 2g-2\\ d=b \,(\mathrm{mod}\ 2)}}\mathcal{M}^{d,c}(\SO_0(2,q)),
\end{equation}
\noi where $\mathcal{M}^{d,c}(\SO(2,q))$ denotes the component in which $\deg(L)=d$ and $sw_2(W)=c$.  

The components where $d=2g-2$ specializes further because in these components:
\begin{enumerate}
\item $L$ has to be isomorphic to $KI$ where $I^2=\mathcal{O}$, and
\item $W$ decomposes as $W=I\oplus W_0$ where $W_0$ is a rank $q-1$ orthogonal bundle with $sw_1(W_0)=I$.
\end{enumerate}

\noi  As shown in \cite{BGG06,CTT}, an $\SO(2,q)$-Higgs bundle with $L=KI$,  $W=I\oplus W_0$ and $\eta=[q_2,\beta]:I\oplus W_0\rightarrow KI$ is defined by a $K^2$-twisted $\SO(1,q-1)$-Higgs bundle $(I,W_0,\beta)$ together with a quadratic differential $q_2$. Denoting $\mathcal{M}^{2g-2,c}(\SO(2,q))$ by $\mathcal{M}_{\mathrm{max}}^c(\SO(2,q))$ it follows  that 
\begin{equation}
\mathcal{M}_{\mathrm{max}}^c(\SO(2,q))=\coprod_{a\in H^1(X,\Z_2)}\mathcal{M}^{a,c}_{K^2}(\SO(1,q-1))\times H^0(K^2)=\coprod_{a}\mathcal{M}^{a,c}_{K^2}(\SO(1,q-1)\times \SO_0(1,1))\end{equation}
\noi where $a=sw_1(I)$.  We thus get
\begin{equation}
\mathcal{M}^{0,0,c}(\SO(2,q))=\coprod_{\substack{0\leq d<2g-2\\ d=0 \,(\mathrm{mod}\ 2)}}
\mathcal{M}^{d,c}(\SO(2,q))\ \sqcup\ \coprod_{a}\mathcal{M}^{a,c}_{K^2}(\SO(1,q-1)\times \SO_0(1,1)).
\end{equation}

\noi The group $\SO(1,q-1)\times\SO_0(1,1)$ is the Cayley partner to $\SO(2,q)$ and the objects in $\mathcal{M}^{a,c}_{K^2}(\SO(1,q-1)\times \SO_0(1,1))$ are the Cayley partners to the Higgs bundles in $\mathcal{M}_{\mathrm{max}}^c(\SO(2,q))$.  As shown in
Section \ref{sec:cayley}
Such Cayley partners are known to emerge in maximal components of $\mathcal{M}(\G)$ whenever $\G$ is the isometry group of a Hermitian symmetric space of tube type. Comparing to Theorem \ref{mainth}, we see that the exotic components in $\mathcal{M}(\SO(p,q))$ are direct generalizations of these Cayley partners to the maximal components $\mathcal{M}_{\mathrm{max}}^c(\SO(2,q))$.

\subsection{Cayley correspondence}

Theorem \ref{mainth} shows not only that additional exotic components exist, but also gives a model which describes them. Indeed, given the model, the result is proved directly by constructing a suitable map from the model to the moduli space $\mathcal{M}(\SO(p,q))$. The model is itself built from moduli spaces (of $K^j$-twisted Higgs bundles), so that in both the domain and target of our map the points represent equivalence classes of objects. one first describes a map between the objects, and then shows that it descends to the appropriate moduli spaces where it defines a homeomorphism onto a connected component.

The map relies in part on a parameterization of the Hitchin components of the moduli spaces $\mathcal{M}(\SO(p-1,p))$.  Viewing theses moduli spaces as subspaces of $\mathcal{M}(\SO(2p-1,\C))$, the parameterization is given by a section of the Hitchin fibration for $\mathcal{M}(\SO(2p-1,\C))$.  As explained in Section \ref{split-slnr}, the fibration is defined by $\SO(2p-1,\C)$-invariant polynomials evaluated on the Higgs field, giving a map to $\bigoplus\limits_{j=1}^{p-1}H^0(K^{2j})$, and admits sections which parameterize connected components of $\mathcal{M}(\SO(p-1,p))\subset\mathcal{M}(\SO(2p-1),\C)$. 

The $\SO(p,p-1)$-Higgs bundles in the image of the section can be taken to be of the form $(\mathcal{K}_{p},\mathcal{K}_{p-1},\sigma(q_2,q_4,\dots,q_{2p-2}))$ where
\begin{equation}\label{Eq cal K notation}
\mathcal{K}_{p}=K^{p-1}\oplus K^{p-3}\oplus\cdots\oplus K^{1-p},
\end{equation}
and $\sigma$ is given by a map
\begin{equation}
	\label{EQ Hitchin section map}
	\sigma:\bigoplus_{j=1}^nH^0(K^{2j})\rightarrow\Hom(\mathcal{K}_{p},\mathcal{K}_{p-1})\otimes K.
\end{equation}

Theorem \ref{mainth} is a consequence of the following (see \cite{A-et-al1,A-et-al2}).

\begin{theorem}\label{cayley-sopq}
 Let $(I,W_0,\eta_p)$ be a $K^p$-twisted $\SO(1,q-p+1)$-Higgs bundle and take differentials $q_{2j}\in  H^0(K^{2j})$ for $j=1,\dots p-1$. Using the notation from \eqref{Eq cal K notation}, consider the $\SO(p,q)$-Higgs bundle 
$(V,W,\eta)$ defined by
\begin{equation}
	\label{EQ etaform}\xymatrix{V=I\otimes \mathcal{K}_p,& W=W_0\oplus I\otimes\mathcal{K}_{p-1}&\text{and}}\ \ \ \ \eta=\mtrx{\mu&\sigma(\ \vec{q}\ )}:W\to V\otimes K,
\end{equation}
where $\mu=\smtrx{\eta_p\\0\\\vdots\\0}:W_0\to \mathcal{K}_p\otimes I\otimes K,$ $\vec{q}=(q_2,\dots,q_{2p-2})$ and $\sigma$ is the map from \eqref{EQ Hitchin section map}.
\begin{enumerate}
\item $(V,W,\eta)$ is polystable if and only if 
 $(I,W_0,\mu)$ is  polystable. The map
\[
 ((I,W_0,\mu); q_2,q_4,\dots,q_{2p-2})\longmapsto (V,W,\eta)\] thus descends to define a map
\[\Psi:\xymatrix{\cM_{K^p}(\SO(1,q-p+1))\times\displaystyle\prod\limits_{i=1}^{p-1} \cM_{K^{2i}}(\SO_0(1,1))\ar[rr]&& \cM(\SO(p,q))}.\]

\item The map $\Psi$ is injective.

\item The image of $\Psi$ is open and closed.
\item The image of $\Psi$ is disjoint from the components $\mathcal{M}_0(\SO(p,q))$.

\end{enumerate}
\end{theorem}

There are several consequences of the results above for  the moduli space of representations $\cR(\SO(p,q))$. Recall that a representation $\rho:\pi_1(S)\to \SO_0(2,1)$ is called {\bf Fuchsian} if it is discrete and faithful and that, since $\SO_0(p,p-1)$ is a split group of adjoint type, there is a unique  principal embedding 
\begin{equation}
\label{eq princ embedd intro}
\iota:\SO_0(2,1)\to\SO_0(p,p-1)~.
\end{equation}

One consequence  is a dichotomy for polystable $\SO(p,q)$-Higgs bundles which translated across the Non-abelian Hodge  correspondence 
leads to the following dichotomy for surface group representations into $\SO(p,q)$.
\begin{theorem} 
Let $S$ be a closed surface of genus $g\geq 2$. For $2<p<q-1$, the moduli space $\cR(\SO(p,q))$ is a disjoint union of two sets,
\begin{equation}
    \label{Eq intro rep dichotomy}
\cR(\SO(p,q))~=~\cR^{cpt}(\SO(p,q))~\sqcup~ \cR^{ex}(\SO(p,q))~,
\end{equation}
where
\begin{itemize}
    \item  $[\rho]\in\cR^{cpt}(\SO(p,q))$ if and only if $\rho$ can be continuously deformed to a compact representation,
    \item  $[\rho]\in\cR^{ex}(\SO(p,q))$ if and only if $\rho$ can be continuously deformed to a representation 
    \begin{equation}\label{eq: model rep intro}
    \rho'=\alpha\oplus (\iota\circ\rho_{\mathrm{Fuch}})\otimes \det(\alpha)~,
    \end{equation}
where $\alpha$ is a representation of $\pi_1(S)$ into the compact group $\OO(q-p+1)$, $\rho_{\mathrm{Fuch}}$ is a Fuchsian representation of $\pi_1(S)$ into $\SO_0(2,1)$, and $\iota$ is the principal embedding from \eqref{eq princ embedd intro}.
 \end{itemize} 
\end{theorem}
\begin{remark}
 For $2<p=q-1$, the above theorem does not hold. Namely, there are exactly $2p(g-1)$ exotic components of $\cR(\SO(p,p+1))$ for which the result fails. 
In  \cite{Col17} Collier conjectures that, with the exception of the Hitchin component,  all representations in these components are Zariski dense. 
\end{remark}

\section{Positivity and Cayley correspondence}

\subsection{Anosov and positive representations}

Anosov representations were introduced by Labourie \cite{labourie} and have many interesting 
geometric and dynamic properties. Important examples of Anosov representations include Fuchsian 
representations, quasi-Fuchsian representations, Hitchin representations into split real groups 
and maximal representations into Lie groups of Hermitian type. We briefly describe the main properties 
of Anosov representations and refer the reader to \cite{labourie,GW12}.

Let $G$ be a semisimple Lie group and $P\subset G$ be a parabolic subgroup. Let $L\subset P$ be the Levi factor
of $P$, it is given by $L=P\cap P^{\opp}$, where $P^{\opp}$ is the opposite parabolic of $P.$ 
The homogeneous space $G/L$ is the unique open $G$ orbit in $G/P\times G/P$, and points $(x,y)\in G/P\times G/P$ in this open orbit are called {\bf transverse}.

\begin{definition}\label{DEF: Anosov rep}
    Let $S$ be  closed surface of genus $g\geq 2$. 
Let $\partial_\infty\pi_1(S)$ be the Gromov boundary of $\pi_1(S)$. Topologically $\partial_\infty\pi_1(S)\cong\R\mathbb{P}^1$. A representation $\rho:\pi_1(S)\to G$ is {\bf  $P$-Anosov} if there exists a unique continuous boundary map $\xi_\rho:\partial_\infty\pi_1(S)\to G/P$
satisfying 
\begin{itemize}
    \item Equivariance: $\xi(\gamma\cdot x)=\rho(\gamma)\cdot\xi(x)$ for all $\gamma\in\pi_1(S)$ and all $x\in\partial_\infty\pi_1(S)$.
    \item Transversality: for all distinct $x,y\in\partial_\infty\pi_1(S)$ the generalized flags $\xi(x)$ and $\xi(y)$ are transverse.
    \item Dynamics preserving: see \cite{labourie,GW12} for the precise notion. 
\end{itemize}
The map $\xi_\rho$ will be called the {\bf $P$-Anosov boundary curve}.
\end{definition}

One important property of Anosov representations is that they define an open subset of the moduli space of representations $\cR(G)$. 
The set of Anosov representations is however not closed. For example, for the group $\PSL(2,\C)$ the set of Anosov representations corresponds to the non-closed set of quasi-Fuchsian representations of $\cR(\PSL(2,\C))$.
The special cases of Hitchin representations and maximal representations define connected components of Anosov representations. Both Hitchin representations and maximal representations satisfy an additional ``positivity'' property which is a closed condition. 
For Hitchin representations this was proved by Labourie \cite{labourie} and Fock--Goncharov \cite{fock-goncharov}, and for maximal 
representations 
by Burger--Iozzi--Wienhard \cite{BIW10}. 
These notions of positivity have recently been unified and generalized by Guichard and Wienhard \cite{GW16}. 

For a parabolic subgroup $P\subset G$, denote the Levi factor of $P$ by $L$ and the unipotent subgroup by $U\subset P$. 
The Lie algebra $\fp$ of $P$ admits an $Ad_{L}$-invariant decomposition $\fp=\fl\oplus\fu$ where $\fl$ and $\fu$ are the Lie 
algebras of $L$ and $U$ respectively. 
Moreover, the unipotent Lie algebra $\fu$ decomposes into irreducible $L$-representation:
\[\fu=\bigoplus\fu_\beta~.\]
Recall that a parabolic subgroup $P$ is determined by fixing a simple restricted root system $\Delta$ of a maximal $\R$-split torus of 
$G$, and choosing a subset $\Theta\subset\Delta$ of simple roots. 
To each simple root $\beta_j\in\Theta$ there is a corresponding irreducible $L$-representation space $\fu_{\beta_j}.$

In \cite{GW16} Guichard and Wienhard introduce the following notions.
A pair $(G,P^\Theta)$ admits a {\bf positive structure} if for all $\beta_j\in\Theta,$ the $L^\Theta$-representation space $\fu_{\beta_j}$ 
has an $L^\Theta_0$-invariant acute convex cone $c_{\beta_j}^\Theta$, where $L^\Theta_0$ denotes the identity component of $L^\Theta$. 

If $(G,P^\Theta)$ admits a positive structure, then exponentiating certain combinations of elements in the $L^\Theta_0$-invariant acute 
convex cones give rise to a semigroup $U^\Theta_{>0}\subset U^\Theta$. The existence of the semigroup $U^\Theta_{>0}$ gives a well defined notion 
of positively oriented triples of pairwise transverse points in $G/P^\Theta.$ This notion allows one to define a {\bf positive Anosov 
representation}. 

If the pair $(G,P^\Theta)$ admits a positive structure, then a $P^\Theta$-Anosov representation $\rho:\pi_1(S)\to\G$ is called 
{\bf positive} if the Anosov boundary curve $\xi:\partial_\infty\pi_1(S)\to G/P^\Theta$ sends positively ordered triples in 
$\partial_\infty\pi_1(S)$ to positive triples in $G/P^\Theta.$

Guichard and Wienhard \cite{GW16} conjecture   that if $(G,P^\Theta)$ admits a notion of positivity, then the 
set $P^\Theta$-positive Anosov representations is an open and  closed subset of $R(G)$.
In particular, the aim of this conjecture is to characterize the connected components of $\cR(G)$ which are not labeled by primary 
topological invariants as being connected components of positive Anosov representations.  This is what defines the  
{\bf higher Teichm\"uller components}. 

When $G$ is a split real form and $\Theta=\Delta$, the corresponding parabolic is a Borel subgroup of $G$. In this case, the connected component of the identity of the Levi factor is $L^\Delta_0\cong(\R^+)^{\rk(\G)}$ and each simple root space $\fu_{\beta_i}$ is one dimensional. The $L^\Delta_0$-invariant acute convex cone in each simple root space $u_{\beta_i}$ is isomorphic to $\R^+.$ The set of $P^\Delta$-positive Anosov 
representations into a split group are exactly Hitchin representations. 
    When $G$ is a Hermitian Lie group of tube type and $P$ is the maximal parabolic associated to the Shilov boundary of the Riemannian 
symmetric space of $G$, the pair $(G,P)$ also admits a notion of positivity \cite{BIW10}. In this case, the space of maximal representations into $G$ are exactly the $P$-positive Anosov representations. In particular, the above conjecture holds in these two cases.

In general, the group $\SO(p,q)$ is not a split group and not a group of Hermitian type. 
Nevertheless, if $p\neq q$, then $\SO(p,q)$ has a parabolic subgroup $P^\Theta$ which admits a positive structure. 
Here $P^\Theta$ is the stabilizer of the partial flag $V_1\subset V_2\subset\cdots \subset V_{p-1},$
where $V_j\subset\R^{p+q}$ is a $j$-plane which is isotropic with respect to a signature $(p,q)$ inner product with $p<q.$ 
Here the subgroup $L^\Theta_{pos}\subset L^\Theta\subset \SO(p,q)$ which preserves the cones $c^\Theta_{\beta_j}$ is isomorphic to 
$L^\Theta_{pos}\cong (\R^+)^l\times \SO(1,q-p+1)$ with $l=p-1$,
We refer the reader to \cite{GW16}  and \cite{Col17} for more details. 

 Using work by Collier \cite{Col17}, the following is shown in  \cite{A-et-al2}. 

\begin{proposition}\label{Prop existence of positive reps}
Let $P^\Theta\subset\SO(p,q)$ be the stabilizer of the partial flag $V_1\subset V_2\subset\cdots \subset V_{p-1},$
where $V_j\subset\R^{p+q}$ is a $j$-plane which is isotropic with respect to a signature $(p,q)$ inner product with $p<q$.  
If $q>p+1$, then each connected component of $\cR^{ex}(\SO(p,q))$ from \eqref{Eq intro rep dichotomy}
contains $P^\Theta$-positive Anosov representations.
\end{proposition}

When $q=p+1$, this was shown in \cite{Col17} for the Collier components mentioned in Section \ref{brian}.

Proposition \ref{Prop existence of positive reps} gives further evidence for Guichard and Wienhard conjecture , and it is thus natural 
to expect that all representations in the connected components from Theorem \ref{cayley-sopq} are positive Anosov 
representations. 
Indeed, this would follow from the conjecture and Proposition \ref{Prop existence of positive reps}. Moreover, if the conjecture 
is true, then the connected components of Theorem \ref{cayley-sopq} correspond exactly to those connected components of 
$\cR(\SO(p,q))$ which contain positive Anosov representations.

\subsection{General Cayley correspondence} 

As we have seen in Section \ref{sec:cayley}, if $G$ is a Hermitian group of 
tube type, subject to a certain topological constraint (always satisfied if $G$ is of adjoint type, for example), there is a bijective correspondence between maximal $G$-Higgs bundles and $K^2$-twisted $H^*$-Higgs bundles, where $H^*$ is a non-compact group determined by the Shilov boundary. 
In fact $H^*$, is the Levi factor of the parabolic subgroup $P$ which defines 
positivity and determines the Shilov boundary as $G/P$.
In particular, Theorem \ref{th:cayley-correspondence} states that  the moduli space of maximal $G$-Higgs bundles is isomorphic to the moduli space of  $K^2$-twisted $H^*$-Higgs bundles.
  
We have also seen in Theorem \ref{cayley-sopq} that for the components of 
the moduli space of 
$\SO(p,q)$-Higgs bundles containing positive representations, for the positive structure defined above, there is also  a generalized Cayley correspondence. Although this seems a bit more involved that the Hermitian case, as explained in
\cite{A-et-al2}, the structure of the Cayley partner is entirely determined 
by the parabolic subgroup of $\SO(p,q)$ defining positivity.    

When $G$ is a split real form the Hitchin components of $\cM(G)$ admit a similar interpretation. The positivity condition in this case is  defined by a minimal parabolic subgroup $P$. The Levi factor is $L=(\R^*)^{\rk(G)}$ and the identity 
component is $L_0=(\R^+)^{\rk(G)}$. Recall from Section \ref{split-slnr} that  the
Hitchin 
base is given
$B(G)=\oplus_{i=1}^r H^0(X,K^{d_i})$, where $r=\rk(G)$. Now, 
the  summand 
$H^0(X,K^{d_i})$ can be interpreted as the moduli space of $K^{d_i}$-twisted 
$\R^+$-Higgs bundles and the Hitchin components are given by
$$
\cM_{K^{d_1}}(\R^+)\times\cdots \times \cM_{K^{d_r}}(\R^+).
$$

From these three situations it seems natural to conjecture that, in general, 
higher  Teichm\"uller components, that is, those consisting of positive 
representations for a certain pair $(G,P)$, are in correspondence  with what 
one can call {\bf Cayley components}, i.e., components of $\cM(G)$ for which 
there is a Cayley correspondence. Moreover,  the structure of the Cayley 
partner  is closely related to the parabolic subgroup defining positivity.
 Work in this direction is being pursued in 
\cite{BCGGO}, building upon results in \cite{BCGT}. This includes   
the exceptional groups for which there is a notion of positivity. As shown 
in \cite{GW16} these  are real forms of $F_4$, $E_6$, $E_7$ and $E_8$ whose 
restricted root system is of type $F_4$.

\appendix
\label{appendix}

\section{Tables}
\label{chap:tables}

We use the following notation for Table \ref{tab:HSS}:

\begin{itemize}
\item $\Delta^\pm_{10}$ are the half-spinor representations of the group
$\Spin(10,\C)$. They are $16$-dimensional.
\item $M$ and $M^*$ are the irreducible $27$-dimensional representations
of $\E_6$, which are dual to each other.
\item $\eta^r$ is the representation $\eta^r: \C^*\to  \C^*$
given by $z\mapsto z^r$.
\end{itemize}

\begin{table}[htbp]
\begin{tabular}{|l|l|l|}
  \hline
Type&${\lieg}$&$\widehat{\lieg}$\\
\hline
AI&$\lie{sl}(n,\R)$&$\sll(n,\R)$\\
\hline
AII&$\lie{su}^*(2n)$&$\sll(n,\R)$\\
\hline
\multirow{2}{*}{AIII}& $\lie{su}(p,q),\ p<q$& $\lie{so}(p,p+1)$\\\cline{2-3}
 &$ \lie{su}(p,p)$ & $\lie{sp}(2p,\R)$\\
\hline
BI&$\lie{so}(2p,2q+1),\ p\leq q$&$\lie{so}(2p,2p+1)$\\
\hline
CI&$\lie{sp}(2n,\R)$&$\lie{sp}(2n,\R)$\\
\hline
\multirow{2}{*}{CII}& $\lie{sp}(p,q)\ p<q$& $\lie{so}(p,p+1)$\\\cline{2-3}
 &$\lie{sp}(p,p)$ & $\lie{sp}(2p,\R)$\\
\hline
BDI&$\lie{so}(p,q)\ p+q=2n, p<q$&$\lie{so}(p,p+1)$\\
\hline
DI&$\lie{so}(p,p)$&$\lie{so}(p,p)$\\
\hline
\multirow{2}{*}{DII}& $\lie{so}^*(4p+2)\ p<q$& $\lie{so}(p,p+1)$\\\cline{2-3}
 &$\lie{so}^*(4p)$ & $\lie{sp}(2p,\R)$\\
\hline
\end{tabular}\qquad
\begin{tabular}{|l|l|l|}
	\hline
	Type&${\lieg}$&$\widehat{\lieg}$\\
	\hline
EI&$\lie{e}_{6(6)}$&$\lie{e}_{6(6)}$\\
\hline
EII&$\lie{e}_{6(2)}$&$\lie{f}_{4(4)}$\\
\hline
EIII&$\lie{e}_{6(-14)}$&$\lie{so}(3,2)$\\
\hline
EIV&$\lie{e}_{6(-26)}$&$\lie{sl}(3,\R)$\\
\hline
EV&$\lie{e}_{7(7)}$&$\lie{e}_{7(7)}$\\
\hline
EVI&$\lie{e}_{7(-5)}$&$\lie{f}_{4(4)}$\\
\hline
EVII&$\lie{e}_{7(-25)}$&$\lie{sp}(6,\R)$\\
\hline
EVIII&$\lie{e}_{8(8)}$&$\lie{e}_{8(8)}$\\
\hline
EIX&$\lie{e}_{8(-24)}$&$\lie{f}_{4(4)}$\\
\hline
FI&$\lie{f}_{4(4)}$&$\lie{f}_{4(4)}$\\
\hline
FII&$\lie{f}_{4(-20)}$&$\lie{sl}(2,\R)$\\
\hline
G&$\lie{g}_{2(2)}$&$\lie{g}_{2(2)}$\\
\hline
\end{tabular}
\vspace{12pt}
\caption{Maximal split subalgebras} \label{tab} 
\end{table}


\begin{center}

\begin{table}[htbp]
\begin{tabular}{|c|c|c|c|}
\hline\raisebox{-8pt}{}
$G$ & $H$ &  $H^\C$ & $\liem^\C=\liem^++ \liem^-$ \\
\hline\hline\raisebox{-8pt}{}
$\SU(p,q)$ & $\SSS(\U(p)\times \U(q))$ & $\SSS(\GL(p,\C)\times
\GL(q,\C))$ &
$\Hom(\C^q,\C^p)+ \Hom(\C^p,\C^q)$ \\
\hline\raisebox{-8pt}{}
$\Sp(2n,\R)$ &  $\U(n)$ & $\GL(n,\C)$ &
$S^2(\C^n) + S^2({\C^n}^*)$ \\
\hline\raisebox{-8pt}{}
$\SO^*(2n)$ &  $\U(n)$  & $\GL(n,\C)$ &
$\Lambda^2(\C^n) + \Lambda^2({\C^n}^*)$ \\
\hline\raisebox{-8pt}{}
$\SO_0(2,n)$ & $\SO(2)\times \SO(n)$ & $\SO(2,\C)\times \SO(n,\C)$ &
 $\Hom(\C^n,\C)+ \Hom(\C,\C^n)$ \\
\hline\raisebox{-8pt}{}
$\E_6^{-14}$ &  $\Spin(10)\times_{\Z_4} \U(1)$  &
$\Spin(10,\C)\times_{\Z_4} \C^*$ &
$\Delta^+_{10}\ot \eta^3 + \Delta^-_{10}\ot \eta^{-3}$ \\
\hline\raisebox{-8pt}{}
$\E_7^{-25}$ & $\E_6^{-78}\times_{\Z_3} \U(1)$  & $\E_6\times_{\Z_3}
\C^*$ &
 $M\ot \eta^2+ M^*\ot \eta^{-2}$ \\

\hline

\end{tabular}
\vspace{12pt}

\caption{Irreducible Hermitian symmetric spaces $G/H$}
    \label{tab:HSS}

\end{table}

\begin{table}[htbp]
\centering
\begin{tabular}{|c|c|c|c|c|c|c|}
\hline\raisebox{-8pt}{}
$G$ & $H$ & $N$ & $\dim \liem$ & $\ell$ & $o(e^{2\pi J})$ & $q_T$\\
\hline\hline\raisebox{-8pt}{}
$\SU(p,q)$ & \footnotesize{$\SSS(\U(p)\times \U(q))$}  & $p+q$ & $2pq$ & \scriptsize{$lcm(p,q)$} & $\frac{p+q}{gcd(p,q)}$ & $\sfrac{1}{2}$ \\
\hline\raisebox{-8pt}{}
$\Sp(2n,\R)$ &  $\U(n)$  & $n+1$ &
 $n(n+1)$ & $n$ & $2$ & $\sfrac{1}{2}$ \\
\hline\raisebox{-8pt}{}
$\SO^*(2n)$ & $\U(n)$  & $2(n-1)$ & $n(n-1)$ & $n$ & $2$ & $1$\\
\hline\raisebox{-8pt}{}
$\SO_0(2,n)$ & $\SO(2)\times \SO(n)$  & $n$ & $2n$ & $1$ & $1$ & $\sfrac{1}{2}$  \\
\hline\raisebox{-8pt}{}
$\E_6^{-14}$ &  \footnotesize{$\Spin(10)\times_{\Z_4} \U(1)$} & $12$ & $32$ & $4$ & $3$ & $\sfrac{1}{2}$\\
\hline\raisebox{-8pt}{}
$\E_7^{-25}$ & $\E_6^{-78}\times_{\Z_3} \U(1)$    & $18$ & $54$ & $3$ & $2$ & $\sfrac{1}{2}$\\
\hline
\end{tabular}
\vspace{12pt}

\caption{Toledo character data for the classical and exceptional groups}
    \label{tab:exponent-tubetype}\label{tab:exponent-nontubetype} \label{tab:exponent-exceptional}
\end{table}

\begin{table}[htbp]
\centering
\begin{tabular}{|c|c|c|c|c|c|c|c|}
\hline\raisebox{-8pt}{}
$G$ & $H$  & $N$ & $\dim \liem$ & $\ell$ & $o(e^{2\pi J})$ & $q_T$\\
\hline\raisebox{-8pt}{}
$\PSU(p,q)$ & \scriptsize{$\mathrm{P}\SSS(\U(p)\times \U(q))$}  & $p+q$ & $2pq$ & \scriptsize{$gcd(p,q)$} & $1$  & $\frac{p+q}{2lcm(p,q)}$ \\
\hline\raisebox{-8pt}{}
\footnotesize{$P\SO^*(2n=4m+2)$} & $\U(n)$ & \small{$2(n-1)$} & \small{$n(n-1)$} & $n$ & $1$ & $2$\\
\hline\raisebox{-8pt}{}
$\E_6^{-14}/\Z_3$ &  \scriptsize{$\Spin(10)\times_{\Z_4} \U(1)$} & $12$ & $32$ & $4$ & $1$ & $\sfrac{3}{2}$\\
\hline
\end{tabular}
\vspace{12pt}
\caption{Toledo character data for adjoint groups of non-tube type.}
    \label{tab:exponent-adjoint}
\end{table}

\newpage

\begin{landscape}

\vspace{12pt}

\vspace{12pt}

\vspace{12pt}

\vspace{12pt}

\renewcommand{\arraystretch}{1.2}

\begin{table}[htbp]
\begin{tabular}{|c|c|c|c|c|c|c|}
\hline\raisebox{-8pt}{}
$G$ & $H$ &  $H^*$ & $H'$ & $\chS=H/H'$ & $\liem'$ &$ \liem'^\C$\\
\hline\hline\raisebox{0pt}{}
$\SU(n,n)$ & $\SSS(\U(n)\times \U(n))$ &
\parbox[t]{110pt}{$\{A \in \GL(n,\C) \st$ \\ \mbox{}\hfill $\det(A)^2 \in \R^+ \}$ \\
}&
\parbox[t]{70pt}{$\{A \in \U(n) \st$ \\ \mbox{}\hfill $\det(A)^2 = 1\}$ \\
}&
$\U(n)$  & $\Herm(n,\C)$  & $\Mat(n,\C)$  \\
\hline\raisebox{-10pt}{}
$\Sp(2n,\R)$ &  $\U(n)$  & $\GL(n,\R)$ & $\OO(n)$ &
$\U(n)/\OO(n)$  & $\Sym(n,\R)$ &  $\Sym(n,\C)$  \\
\hline\raisebox{-20pt}{}
\parbox[t]{40pt}{$\SO^*(2n)$\\ \mbox{}\hfill  $n=2m$} &  $\U(n)$  & $\U^*(n)$ & $\Sp(n)$ &
$\U(n)/\Sp(n)$  & $\Herm(m,\mathbb{H})$ &  $\Skew(n,\C)$  \\
\hline\raisebox{-20pt}{}
$\SO_0(2,n)$ & $\SO(2)\times \SO(n)$ & $\SO_0(1,1)\times \SO(1,n-1)$
& $\OO(n-1)$ &
\raisebox{-5pt}{$\dfrac{\U(1)\times S^{n-1}}{\Z_2}$}  & $\R\times \R^{n-1}$ &  $\C\times \C^{n-1}$  \\
\hline\raisebox{-20pt}{}
$\E_7^{-25}$ & $\E_6^{-78}\ts_{\Z_3} \U(1)$ & $\E_6^{-26}\ltimes \R^*$ & $\F_4\times \Z_2$ &
\raisebox{-5pt}{$\dfrac{\E_6^{-78}\cdot \U(1)}{\F_4}$}  & $\Herm(3,\Oc)$ &  $\Herm(3,\Oc)\otimes \C$  \\
\hline
\end{tabular}
\vspace{10pt}
\caption{Irreducible Hermitian symmetric spaces
$G/H$ of tube type}
    \label{tab:tube}
\end{table}


\end{landscape}

\end{center}

\clearpage

\newcommand{\etalchar}[1]{$^{#1}$}


\begin{thebibliography}{AAAAA00}

\bibitem[ABC+18]{A-et-al1} M. Aparicio-Arroyo, S.B. Bradlow, B. Collier, 
O. Garc\'{\i}a-Prada, P.B. Gothen and A. Oliveira,
Exotic components of $\SO(p,q)$ surface group representations and their 
Higgs bundle avatars, {\em C.R. Math. Acad. Sci. Paris} {\bf 356} (2018) 
666--673.

\bibitem[ABC+18b]{A-et-al2} 
\bysame, $\SO(p,q)$-Higgs bundles and higher Teichm\"uller components,
\texttt{arXiv:1802.08093}, 2018.

\bibitem[AB83]{AB83}
M.~F. Atiyah and R.~Bott \newblock The {Y}ang-{M}ills equations over {R}iemann surfaces.
\newblock {\em Philos. Trans. Roy. Soc. London Ser. A}, {\bf 308}
(1983) 523--615.

\bibitem[BCGT]{BCGT} O. Biquard, B. Collier,  O. Garc{\'{\i}}a-Prada and
D. Toledo, Arakelov--Milnor inequalities and maximal variations of Hodge 
structure. In preparation.

\bibitem[BGM15]{BGM15} O. Biquard,  O. Garc{\'{\i}}a-Prada and I. Mundet 
i Riera, Parabolic Higgs bundles and representations of the fundamental group 
of a puntured surface into a real group, \texttt{arXiv:1510.042207}.

\bibitem[BGR17]{BGR17} O. Biquard,  O. Garc{\'{\i}}a-Prada and R. Rubio,
Higgs bundles, the Toledo invariant and the Cayley correspondence, 
{\em J. of Topology} {\bf 10} (2017) 795--826.

\bibitem[BCGGO]{BCGGO} 
S.B. Bradlow, B. Collier,  O. Garc\'{\i}a-Prada, P.B. Gothen and A. Oliveira,
Higher Teichm\"uller spaces, positivity and  Cayley correspondence. 
In preparation.

\bibitem[BGG03]{BGG03}
S.B. Bradlow, O. Garc{\'{\i}}a-Prada, and P.B. Gothen.
\newblock Surface group representations and {${\rm U}(p,q)$}-{H}iggs bundles.
\newblock {\em J. Differential Geom.} {\bf 64} (2003) 111--170.

\bibitem[BGG06]{BGG06}
\bysame,
\newblock Maximal surface group representations in isometry groups of classical
 {H}ermitian symmetric spaces.
\newblock {\em Geom. Dedicata}, {\bf 122} (2006) 185--213.

\bibitem[BGG13]{BGG13}
\bysame,
Deformations of maximal 
representations in $\Sp(4,\R)$, {\em Q. J. Math.} \textbf{63} (2012) 795--843.

\bibitem[BGG15]{BGG15}
\bysame,
Higgs bundles for the non-compact dual of the special orthogonal group
{\em  Geom. Dedicata} {\bf 175} (2015) 1--48.

\bibitem[BGM03]{BGM03}
S.~B. Bradlow, O.~Garc{\'{\i}}a-Prada, and I.~Mundet~i Riera, Relative
{H}itchin-{K}obayashi correspondences for principal pairs,
{\em Quarterly J. Math.},   \textbf{54} (2003) 171--208.

\bibitem[BILW05]{burger-iozzi-labourie-wienhard:2005}
M.~Burger, A.~Iozzi, F.~Labourie, and A.~Wienhard, Maximal
representations of surface groups: symplectic {A}nosov structures, {\em Pure
Appl. Math. Q.} \textbf{1} (2005) 543--590.

\bibitem[BIW10]{BIW10}
M. Burger, A. Iozzi, and A. Wienhard.
\newblock Surface group representations with maximal {T}oledo invariant.
\newblock {\em Ann. of Math. (2)}  {\bf 172} (2010) 517--566.

\bibitem[BIW14]{BIW14}
M.~{Burger}, A.~{Iozzi}, and A.~{Wienhard}, Higher Teichm\"{u}ller Spaces: 
from $\SL(2,\R)$ to other Lie  groups. In {\em Handbook of Teichm\"uller 
theory,  Vol. IV}, volume 19 of IRMA Lect. Math. Theor. Phys., pages
539--618. EMS, Z\"urich 2014

\bibitem[ChG93]{ChG93} S. Choi and W.M. Goldman, Convex real projective structures
on closed surfaces are closed {\em Proc. Amer. Math. Soc.}{\bf 118} (1993) 657--661.  

\bibitem[CO03]{CO03}
J.-L. Clerc and B. Orsted, The Gromov norm of the Kaehler class and the
Maslov index {\em  Asian J. Math.}, {\bf 7} (2003) 269--295.

\bibitem[Col17]{Col17} B. Collier, $\SO(n,n+1)$-surface group representations 
and their Higgs bundles, \texttt{arXiv:1710.01287}, 2017.

\bibitem[CTT17]{CTT} B. Collier, N. Tholozan, J. Toulisse, 
The geometry of maximal 
representations of surface groups into $\SO(2,n)$, \texttt{arXiv:1702.08799}, 2017.

\bibitem[Cor88]{corlette}
K.~Corlette, Flat ${G}$-bundles with canonical metrics, \textsl{J. Differential
  Geom.}, \textbf{28} (1988) 361--382.

\bibitem[DG01]{DG01}  R.Y. Donagi and D. Gaitsgory,
 The gerbe of Higgs bundles,  {\em Transform.
Groups}, {\bf 7} (2001) 109--153.

\bibitem[Don87]{donaldson}
S.K. Donaldson, Twisted harmonic maps and the self-duality equations.
{\em Proc. London Math. Soc.} (3) \textbf{55} (1987) 127--131.

\bibitem[DT87]{DT87}
A. Domic and  D. Toledo, The {G}romov norm of the {K}aehler class of symmetric
domains. {\em Math. Ann.} \textbf{276} (1987) 425--432.

\bibitem[Dup78]{dupont}
J. L. Dupont, {\em Bounds for characteristic numbers of 
flat bundles}, Springer LNM, {\bf 763} (1978).

\bibitem[FK94]{FK94}
J. Faraut and A. Kor{\'a}nyi.
\newblock {\em Analysis on symmetric cones}.
\newblock Oxford Mathematical Monographs. The Clarendon Press Oxford University
  Press, New York, 1994.
\newblock Oxford Science Publications.

\bibitem[FG06]{fock-goncharov} V.V. Fock and A.B. Goncharov, Moduli spaces of 
local systems and higher Teichm\"uller theory,
{\em Inst. Hautes {\'E}tudes Sci. Publ. Math.} \textbf{103} (2006) 1--211.
 
\bibitem[GGM09]{GGM09}
O.~{Garc\'{\i}a-Prada}, P.~B. {Gothen}, and I. Mundet i {Riera}.
\newblock {The Hitchin-Kobayashi correspondence, Higgs pairs and surface group
  representations},  \texttt{arXiv:0909.4487}.

\bibitem[GGM13]{GGM13} 
\bysame,
Higgs bundles and surface group representations in the real
symplectic group, {\em Journal of Topology}, {\bf 6} (2013) 64--118.

\bibitem[GM04]{GM04}
O.~Garc{\'\i}a-Prada and I. Mundet i Riera, Representations of the fundamental
group of a closed oriented surface in $\Sp(4,\R)$, 
{\em Topology}, \textbf{43} (2004) 831--855.

\bibitem[GO14]{garcia-prada-oliveira}
O. Garc\'{\i}a-Prada
and  A. Oliveira, 
Connectedness of Higgs bundle moduli for complex reductive Lie groups,
{\em Asian Journal of Mathematics}, {\bf 21} (2017) 791--810.

\bibitem[GP19]{GP19}
O. Garc\'{\i}a-Prada and A. Pe\'on-Nieto, Higgs pairs, abelian gerbes and 
cameral data, preprint 2019.

\bibitem[GPR18]{GPR18}
O. Garc\'{\i}a-Prada, A. Pe\'on-Nieto and S. Ramanan, Higgs bundles for real 
groups and the Hitchin--Kostant--Rallis section {\em Transactions of the AMS} 
{\bf 370} (2018) 2907--2953.

\bibitem[GR18]{garcia-prada-ramanan}
 O. Garc\'{\i}a-Prada and S. Ramanan,
Involutions of rank 2  Higgs bundle moduli spaces. In 
{\em Geometry and Physics:  A Festschrift in honour of Nigel Hitchin}, 
Oxford University Press, 2018.

\bibitem[GW08]{guichard-wienhard:2008}
O.~Guichard and A.~Wienhard, Convex foliated projective structures and
the {H}itchin component for $\PSL(4,\R)$, {\em Duke Math. J.}
\textbf{144} (2008) 381--445.

\bibitem[GW12]{GW12}
O.~Guichard and A.~Wienhard, Anosov representations: domains of discontinuity 
and applications  {\em Invent. Math.}
\textbf{190} (2012) 357--438.

\bibitem[GW16]{GW16}
O.~Guichard and A.~Wienhard, Positivity and higher Teichm\"uller theory, 
{\em Proceedings of the 7th European Congress of Mathematics}, 2016.

\bibitem[Gol88]{goldman}
W.M. Goldman,  Topological components of spaces of representations.
{\em Invent. Math.} \textbf{93} (1988) 557--607.

\bibitem[Got01]{gothen}
P. B. Gothen,  Components of spaces of representations and stable triples,
{\em Topology} \textbf{40} (2001) 823--850.

\bibitem[Hel01]{helgason}
S. Helgason,
\newblock {\em Differential geometry, {L}ie groups, and symmetric spaces},
  volume~34 of {\em Graduate Studies in Mathematics}.
\newblock American Mathematical Society, Providence, RI, 2001.
\newblock Corrected reprint of the 1978 original.

\bibitem[Her91]{hernandez}
L. Hern\'andez, Maximal representations of surface groups in bounded symmetric
domains {\em  Trans. Amer. Math. Soc.}, {\bf 324} (1991) 405--420.

\bibitem[Hit87]{hitchin1987}
N.~J. Hitchin, The self-duality equations on a {R}iemann surface,
{\em Proc. London Math. Soc.}, \textbf{55} (1987) 59--126.

\bibitem[Hit87b]{hitchin:duke}
\bysame,
Stable bundles and integrable systems, {\em Duke Math. J.} {\bf 54} 
(1987) 91--114.

\bibitem[Hit92]{hitchin92}
\bysame,
\newblock Lie groups and {T}eichm\"uller space.
\newblock {\em Topology}, {\bf 31}  (1992) 449--473.

\bibitem[HS14]{HS14} N.~J. Hitchin and L. Schaposnik, 
Nonabelianization of Higgs bundles, 
{\em J. Differential Geom.} {\bf 97}, (2014) 79--89.

\bibitem[Kna96]{knapp}
A. W. Knapp, {\em Lie Groups beyond an Introduction}, first ed.,
Progress in Mathematics, vol 140, Birkh\"auser Boston Inc., Boston,
MA, 1996.

\bibitem[KP14]{KP14} I. Kim and P. Pansu, Density of Zariski density for 
surface groups, {\em Duke Math. J.} {\bf 163} (2014) 1737--1794.

\bibitem[Kob87]{kobayashi}
S. Kobayashi.
\newblock {\em Differential geometry of complex vector bundles}, volume~15 of
  {\em Publications of the Mathematical Society of Japan}.
\newblock Princeton University Press, Princeton, NJ, 1987.
\newblock Kan{\^o} Memorial Lectures, 5.

\bibitem[KV79]{KV79}
A. Kor{\'a}nyi and S. V{\'a}gi.
\newblock Rational inner functions on bounded symmetric domains.
\newblock {\em Trans. Amer. Math. Soc.}, \textbf{254} (1979) 179--193.


\bibitem[KR71]{KR71} B. Kostant and S.  Rallis, Orbits and representations 
associated with symmetric spaces,
{\em Amer. J. Math.} {\bf 93} (1971) 753--809.


\bibitem[Lab06]{labourie}
F. Labourie, Anosov flows, surface groups and curves in projective space,
{\em Invent. Math.} \textbf{165} (2006) 51--114.

\bibitem[Li93]{li}
J. Li, The Space of Surface Group Representations, 
\textsl{Manuscripta Math.} \textbf{78} (1993), 223--243.

\bibitem[Mil58]{milnor}
J.~W. Milnor, On the existence of a connection with curvature zero,
{\em Commm. Math. Helv.} {\bf 32} (1958) 215--223.

\bibitem[Nit91]{nitsure}
 N. Nitsure, Moduli space of semistable pairs on a
curve, {\em Proc. London Math. Soc.} {\bf 62} (1991) 275--300.


\bibitem[Peo13]{Peo13} A. Pe\'on-Nieto, {\em Higgs bundles, real forms and the Hitchin fibration}, Ph.D. Thesis, UAM-ICMAT, Madrid 2013.


\bibitem[Ram96]{ramanathan96} 
A.~Ramanathan.
 Moduli for principal
  bundles over algebraic curves: I and II, {\em Proc. Indian
Acad. Sci. Math. Sci.}  \textbf{106} (1996) 301--328 and 421--449.


\bibitem[Sch08]{schmitt08}
A. H.~W. Schmitt.
\newblock {\em Geometric invariant theory and decorated principal bundles}.
\newblock Zurich Lectures in Advanced Mathematics. European Mathematical
  Society (EMS), Z\"urich, 2008.

\bibitem[Sim88]{Sim88}
C.~T. Simpson, Constructing variations of {H}odge structure using
  {Y}ang--{M}ills theory and applications to uniformization, {\em J. Amer. Math.
  Soc.}  \textbf{1} (1988) 867--918.

\bibitem[Sim90]{Sim90}
C.~T. Simpson, 
Harmonic bundles of non-compact curves
{\em J. Amer. Math. Soc.}  \textbf{3} (1990) 713--770.

\bibitem[Sim92]{Sim92}
\bysame, 
Higgs bundles and local systems, {\em Inst. Hautes {\'E}tudes Sci.
  Publ. Math.} \textbf{75} (1992), 5--95.

\bibitem[Sim94]{simpson94}
\bysame, 
{Moduli of representations of the fundamental group of a smooth
  projective variety {I}}, {\em Publ. Math., Inst. Hautes \'Etud. Sci.}
 \textbf{79}   (1994) 47--129.

\bibitem[Sim95]{simpson95}
\bysame, {Moduli of representations of the fundamental group of a smooth
 projective variety {II}}, {\em Publ. Math., Inst. Hautes \'Etud. Sci.}
 \textbf{80} (1995) 5--79.

\bibitem[Sch13]{Sch13}
L.P. Schaposnik, {\em Spectral data for $G$-Higgs bundles}, DPhil thesis, 
University of Oxford, 2013.

\bibitem[Tol89]{toledo}
D. Toledo.
\newblock Representations of surface groups in complex hyperbolic space.
\newblock {\em J. Differential Geom.}, \textbf{29} (1989) 125--133.

\bibitem[Wien18]{wienhard}
A. Wienhard, An invitation to higher Teichm\"uller theory, 
{\em Proc. Int. Cong. of Math.} 
Rio de Janeiro, Vol. 1 (2018), 1007--1034.

\bibitem[Woo71]{wood}
J. W. Wood, Bundles with totally disconnected structure group,
{\em  Comment. Math. Helv.} {\bf 46} (1971) 257--273.

\end{thebibliography}
\end{document}